\documentclass[twoside,11pt]{article}

\usepackage[preprint]{jmlr2e}

\usepackage{jmlr2e}

\usepackage{hyperref} 

\usepackage{times,tabularx,subfigure}
\usepackage{amsfonts,amsmath,mathabx,caption}
\usepackage{enumerate,wrapfig}
\usepackage[shortlabels]{enumitem}
\usepackage{multicol}
\usepackage{empheq}
\usepackage{color,makecell,booktabs}
\usepackage{algorithm,algorithmic}

\newcommand{\qed}{\hfill \ensuremath{\Box}}

\allowdisplaybreaks
\let\underbrace\LaTeXunderbrace

\newcommand{\vct}[1]{\boldsymbol{#1}}

\newcommand{\vtheta}{\vct{\theta}}

\usepackage{tikz}

\newcommand{\Acal}{{\cal A}}
\newcommand{\Bcal}{{\cal B}}

\newcommand{\Fcal}{{\cal F}}

\newcommand{\Ocal}{{\cal O}}
\newcommand{\Pcal}{{\cal P}}

\newcommand{\Scal}{{\cal S}}

\newcommand{\Ucal}{{\cal U}}

\newcommand{\Xcal}{{\cal X}}

\newcommand{\ybf}{{\bf y}}

\newtheorem{lem}{Lemma}
\newtheorem{thm}{Theorem}

\newtheorem{assump}{Assumption}

\usepackage{lastpage}

\ShortHeadings{Accelerated Multi-Time-Scale Stochastic Approximation}{Zeng and Doan}
\firstpageno{1}

\begin{document}

\title{Accelerated Multi-Time-Scale Stochastic Approximation: Optimal Complexity and Applications in Reinforcement Learning and Multi-Agent Games}

\author{
\name Sihan Zeng \email szeng2017@gmail.com \\
       \addr JPMorgan AI Research\\
       Palo Alto, CA, United States
       \AND
       \name Thinh T. Doan \email thinhdoan@utexas.edu \\
       \addr Department of Aerospae Engineering \& Engineering Mechanics\\
       University of Texas at Austin\\
       Austin, TX, United States
}

\editor{My editor}
 
\maketitle

\begin{abstract}
Multi-time-scale stochastic approximation is an iterative algorithm for finding the fixed point of a set of $N$ coupled operators given their noisy samples. It has been observed that due to the coupling between the decision variables and noisy samples of the operators, the performance of this method decays as $N$ increases. In this work, we develop a new accelerated variant of multi-time-scale stochastic approximation, which significantly improves the convergence rates of its standard counterpart.   
Our key idea is to introduce auxiliary variables to dynamically estimate the operators from their samples, which are then used to update the decision variables. These auxiliary variables help not only to control the variance of the operator estimates but also to decouple the sampling noise and the decision variables. This allows us to select more aggressive step sizes to achieve an optimal convergence rate. Specifically, under a strong monotonicity condition, we show that for any value of $N$ the $t^{\text{th}}$ iterate of the proposed algorithm converges to the desired solution at a rate $\widetilde{\Ocal}(1/t)$ when the operator samples are generated from a single from Markov process trajectory.  
%

A second contribution of this work is to demonstrate that the objective of a range of problems in reinforcement learning and multi-agent games can be expressed as a system of fixed-point equations. As such, the proposed approach can be used to design new learning algorithms for solving these problems. 
We illustrate this observation with numerical simulations in a multi-agent game and show the advantage of the proposed method over the standard multi-time-scale stochastic approximation algorithm.
\end{abstract}


\section{Introduction}

The objective of a range of problems in reinforcement learning and games can be expressed as a coupled system of $N$ equations, with each equation defined through a nonlinear operator that can be queried with noise. 
This work proposes a single-loop stochastic approximation (SA) algorithm which solves such a system using a (near) optimal number of queries. 
Given a compact statistical sample space $\Xcal$ and operators $\{F_i:\mathbb{R}^{d_1}\times\cdots\times\mathbb{R}^{d_N}\times\Xcal\rightarrow\mathbb{R}^{d_i}\}_{i=1,\cdots,N}$, our aim is to find a solution tuple $\vtheta^{\star}=(\theta_1^{\star},\cdots,\theta_N^{\star})\in\mathbb{R}^{d_1}\times\cdots\times \mathbb{R}^{d_N}$ such that
\begin{align}
\left\{\;
\begin{array}{c}
     \mathbb{E}_{X\sim\mu_{\vtheta^{\star}}}[F_1(\theta_1^{\star},\theta_2^{\star},\cdots,\theta_N^{\star},X)]=0, \\
     \mathbb{E}_{X\sim \mu_{\vtheta^{\star}}}[F_2(\theta_1^{\star},\theta_2^{\star},\cdots,\theta_N^{\star},X)]=0, \\
     \vdots \\
     \mathbb{E}_{X\sim \mu_{\vtheta^{\star}}}[F_N(\theta_1^{\star},\theta_2^{\star},\cdots,\theta_N^{\star},X)]=0,
\end{array}
\right.
\label{eq:obj}
\end{align}
where $\mu_{\vtheta}$ denotes a distribution over $\Xcal$ (possibly) as a function of $\vtheta$.
We do not assume having direct access to $\mu_{\vtheta}$. Rather, our sampling oracle is that we may query an operator $P_{\vtheta}:\Xcal\rightarrow\Delta_{\Xcal}$, which induces $\mu_{\vtheta}$ as the stationary distribution. More precisely, for any $\vtheta$ and $X$, we can draw a sample $X'\sim P_{\vtheta}(X)$.
We consider this relaxed sampling oracle as it realistically captures how samples are generated in reinforcement learning (RL).
We use $\Omega=\mathbb{R}^{d_1}\times\cdots\times \mathbb{R}^{d_N}$ to denote the space of the decision variable $\vtheta$ and denote $\widebar{F}_i(\theta_1,\cdots,\theta_N)=\mathbb{E}_{X\sim\mu_{\vtheta}}[F_i(\theta_1,\cdots,\theta_N,X)]$ for any $\vtheta\in\Omega$.\looseness=-1

It may appear that the system \eqref{eq:obj} is symmetric and
each equation is identically coupled with any other equation.
This is not true in the context of RL and games, where the problems usually have a specific hierarchical structure. The structure dictates the order in which the system can be most naturally solved and is formalized later through Assumption~\ref{assump:stronglymonotone}.

Exploiting the hierarchical structure, the standard multi-time-scale stochastic approximation algorithm for solving \eqref{eq:obj}, thereafter referred to as \textbf{\texttt{MSA}}, maintains parameters $\theta_1^{[t]},\cdots,\theta_N^{[t]}$ as an estimate of $(\theta_1^{\star},\cdots,\theta_N^{\star})$ and simultaneously updates them, each using a unique step size. 
Given (possibly biased) stochastic samples $\{X^{[t]}\}$, \texttt{MSA} iteratively carries out
\begin{align}
    \theta_i^{[t+1]} = \theta_i^{[t]} - \alpha_i^{[t]} F_i(\theta_1^{[t]},\cdots,\theta_N^{[t]},X^{[t]}),\quad\forall i=1,\cdots,N.\label{alg:update_standardSA}
\end{align}
With properly selected step sizes, this simple method is guaranteed to converge to the solution of \eqref{eq:obj} under strong monotonicity and Lipschitz continuity assumptions on $\widebar{F}_i$ to be introduced later \citep{hong2023two,zeng2024two,shen2022single,han2024finite}. However, its time and sample complexity degrades as $N$ increases and is sub-optimal as known so far. 
When $N=1$, \texttt{MSA} reduces to the classic stochastic approximation algorithm \citep{robbins1951stochastic} and has a complexity of $\Ocal(1/t)$ when $X^{[t]}$ are i.i.d. samples from a fixed distribution (and $\widetilde{\Ocal}(1/t)$ if the samples $X^{[t]}$ are consecutively generated from a Markov chain, where $\widetilde{\Ocal}$ hides any logarithm factor). More precisely, the last iterate of \eqref{alg:update_standardSA} converges to an $\widetilde{\Ocal}(1/t)$ neighborhood around $\vtheta^{\star}$ in $t$ iterations, matching the worst-case lower bound \citep{chen2022finite}.
When $N=2$, the best-known complexity, as established in \citet{hong2023two,zeng2024two}, becomes $\widetilde{\Ocal}(1/t^{2/3})$, inferior to that under $N=1$ by a factor of $t^{1/3}$. The authors in \citet{hong2023two,zeng2024two} only analyze the case of two equations, but we show in Section~\ref{sec:standardSA} that their techniques can be extended to establish a complexity of $\widetilde{\Ocal}(1/t^{1/2})$ for $N=3$ and further worse bounds for $N=4,5,\cdots$. This apparent gap in the performance of \texttt{MSA} for different $N$ makes us question whether the problem becomes fundamentally harder as $N$ increases or whether any algorithmic/analytical improvement can be made. In this work, we answer the question in the positive direction by proposing an accelerated SA algorithm which achieves the optimal complexity $\widetilde{\Ocal}(1/t)$ across all $N$ while maintaining the favorable properties of \eqref{alg:update_standardSA} such as simplicity and single loop structure. The same algorithm has been introduced earlier in our prior works \citep{doan2024fast,zeng2024fast} but we make important contributions over them as discussed in details in Section~\ref{sec:literature}.
We now summarize our main contributions.

\subsection{Main Contributions}

\begin{wraptable}{r}{0.5\textwidth}
\begin{tabular}{ccc}
\toprule
\makecell{Number of\\ Equations}  & \makecell{Standard \\
Multi-Time-Scale\\ SA Complexity} & \makecell{\texttt{A-MSA}\\  Complexity}\\
\midrule
Two & $\widetilde{\Ocal}(\epsilon^{-1.5})$ & $\widetilde{\Ocal}(\epsilon^{-1})$\\
Three & $\widetilde{\Ocal}(\epsilon^{-2})$ & $\widetilde{\Ocal}(\epsilon^{-1})$\\
Four & $\widetilde{\Ocal}(\epsilon^{-2.5})$ & $\widetilde{\Ocal}(\epsilon^{-1})$ \\
$N$ & $\widetilde{\Ocal}(\epsilon^{-(N+1)/2})$ & $\widetilde{\Ocal}(\epsilon^{-1})$\\
\bottomrule
\bottomrule
\end{tabular}
\caption{Sample complexity comparison as number of equations increases (under strong monotone operators and Markovian samples).}
\label{table:comparison}
\end{wraptable}

The key contribution of this work is to propose a novel Accelerated Multi-Time-Scale Stochastic Approximation (\textbf{\texttt{A-MSA}}) algorithm that solves a coupled system of 
$N$ equations with the optimal convergence rate $\widetilde{\Ocal}(1/t)$ for strong monotone operators and when the samples are generated from a single trajectory of Markov processes characterized by $P_{\vtheta}$. In the case $N=2$, this rate improves over $\widetilde{\Ocal}(1/t^{2/3})$, the best-known complexity of the \texttt{MSA} algorithm derived in \citet{doan2022nonlinear,hong2023two}. When $N\geq 3$, the complexity of \texttt{MSA} deteriorates and is further inferior to our result (see Table.~\ref{table:comparison}).
Our innovation is built on the observation that the main hindrance to the fast convergence of \eqref{alg:update_standardSA} lies in the direct coupling between the iterates $\theta_i^{[t]}$ across $i$, which means that any error/noise in $\theta_i^{[t]}$ for some $i$ immediately propagates to $\theta_j^{[t+1]}$ for all $j$ in the next iteration. Such a coupling effect prevents us from choosing aggressive step sizes necessary for achieving the optimal complexity. Our solution to this challenge is to update $\theta_i^{[t]}$ in the direction of a properly averaged version of operator $F_i$ (as opposed to the single sample used in \eqref{alg:update_standardSA}). The seemingly small modification maintains the simplicity of \eqref{alg:update_standardSA} while eliminating the direct coupling between $\theta_i^{[t]}$. In Section~\ref{sec:comparison}, we compare the analysis of \texttt{MSA} and \texttt{A-MSA} and pinpoint the convergence bottleneck in \texttt{MSA} which we overcome through an algorithmic modification.

A second contribution of our work is to show that coupled equations in \eqref{eq:obj} abstract many popular problems in games and RL, including gradient-based temporal difference learning with momentum, distributionally robust RL, policy optimization in a two-player zero-sum Markov game, and policy optimization in mean-field games. As a consequence, the proposed \texttt{A-MSA} method can be used to design new algorithms for solving these problems. Some of these problems do not satisfy the strong monotonicity condition but are specially structured in their own ways. We show how we can adapt the analysis for strongly monotone operators to the structure in each respective problem, which can potentially establish the state-of-the-art/previously unknown convergence rates.


\begin{table}[!ht]
\centering
\setlength{\tabcolsep}{1pt}
\begin{tabular}{ccccc}
        \toprule
        & \makecell{Applicable to\\ $>2$ Equations}  & \makecell{Sample \\ Noise} & \makecell{Smoothness \\ Assumption} & \makecell{Sample \\ Complexity}\\
        \midrule
        \citet{deb2021n} & Yes & Bounded Martingale & Not Required & -\\
        \citet{hong2023two} & No & i.i.d. & Not Required & $\Ocal(\epsilon^{-1.5})$\\
        \citet{zeng2024two} & No & Time-Varying Markovian & Not Required & $\widetilde{\Ocal}(\epsilon^{-1.5})$\\
        \citet{shen2022single} & Yes & Bounded Martingale & Required & $\Ocal(\epsilon^{-1})$\\
        \citet{han2024finite} & No & Bounded Martingale & Relaxed Version & $\Ocal(\epsilon^{-1})$\\
        \citet{doan2024fast} & No & i.i.d. & Not Required & $\Ocal(\epsilon^{-1})$\\
        \citet{zeng2024fast} & No & i.i.d. & Not Required & $\Ocal(\epsilon^{-1})$\\
        \textbf{This Work} & \textbf{Yes} & \textbf{Time-Varying Markovian} & \textbf{Not Required} & \textbf{$\widetilde{\Ocal}(\epsilon^{-1})$}\\
        \bottomrule
        \bottomrule
        \end{tabular}
\caption{Existing algorithms on nonlinear two/multi-time-scale stochastic approximation and their assumptions and complexity.}
\vspace{-10pt}
\end{table}

\subsection{Related Work}\label{sec:literature}

This paper closely relates to the existing works on singe- and multi-time-scale SA (especially under Markovian samples) and those that analyze the finite-time and sample complexity of reinforcement learning algorithms through the lens of SA. We discuss the most relevant works in these domains and highlight our novelty.

\noindent\textbf{Stochastic Approximation.} The SA algorithm is a classic method for solving a single fixed-point equation given noisy operator samples \citep{robbins1951stochastic}. The asymptotic convergence of SA is well understood and usually derived with the Poisson equation method \citep{benveniste2012adaptive,li2023online} or by analyzing a related ODE \citep{Meerkov1972a,borkar2008stochastic}. 
Under i.i.d. samples, bounded Martingale difference noise, or dependent samples from a Markov chain with bounded mixing time, the finite-time analysis measured by squared errors in expectation is established for linear \citep{lakshminarayanan2017linear,srikant2019finite} and strongly monotone operators \citep{moulines2011non,chen2022finite}. 
Linear SA abstracts the temporal difference (TD) learning algorithm in RL. Various works \citep{bhandari2018finite,mitra2024simple} study the finite-time convergence of TD learning leveraging the problem structure, but their analysis may be extended to linear SA. Recently, non-asymptotic central limit theorems have also been proved for TD learning \citep{srikant2024rates} and Q learning \citep{li2023statistical}, which show the error convergence to a zero-mean Gaussian random variable in distribution. Such results may be extended to general linear/non-linear SA as well.

Some of the aforementioned works employ Polyak-Ruppert averaging to improve variance control \citep{srikant2024rates,haque2023tight}. We note that the averaging step we carry out in the proposed method is distinct from Polyak-Ruppert averaging -- Polyak-Ruppert averaging does not modify the algorithm itself but returns the solution as an weighted-averaged version of all history iterates, while our algorithm uses an averaged sequence in the updates.

\noindent\textbf{Multi-Time-Scale Stochastic Approximation.} When $N=2$, the MSA algorithm reduces to the standard two-time-scale SA, which is known to converge asymptotically under linear \citep{borkar1997stochastic,borkar2008stochastic,konda2004convergence} and non-linear operators \citep{mokkadem2006convergence,fort2015central}.
In the linear setting, the finite-time convergence measured by mean squared errors is first established with a sub-optimal rate $\widetilde{O}(1/t^{2/3})$ in \citet{gupta2019finite,doan2019linear} and improved to the optimal $\widetilde{O}(1/t)$ in \citet{kaledin2020finite,haque2023tight}. In addition, convergence in probability is studied in \cite{dalal2018finite}. In the non-linear setting under the strong monotonicity assumption, the best-known rate of two-time-scale SA remains $\widetilde{O}(1/t^{2/3})$ in general \citep{doan2022nonlinear, hong2023two,zeng2024two}. Under an additional smoothness assumption, \citet{shen2022single,han2024finite} show that the optimal rate can be recovered. Our paper does not require this assumption as it may not hold true in practical applications.

The two-time-scale SA framework has wide applications in control, RL, and games, e.g., gradient-based TD learning \citep{sutton2009fast,wang2021non},  actor-critic methods \citep{konda1999actora,konda1999actorb,wu2020finite}, policy gradient descent-ascent for two-player games \citep{daskalakis2020independent}, and decentralized Q-learning \citep{sayin2021decentralized}. 
By providing an analysis for the general framework, one can apply it to each specific problem and immediately deduce the performance guarantees, without having to tailor the analysis on a case basis.

The study of SA for solving three or more equations remains limited but deserves increased attention due to its ability to model rich classes of problems, including gradient-based TD learning with momentum \citep{deb2022gradient}, distributionally robust RL \citep{liang2024single}, and policy optimization in games; the applications will be discussed in Section~\ref{sec:applications}. The existing works on multi-time-scale SA include \citet{deb2021n} which only provides an asymptotic analysis and \citet{shen2022single} which relies on a restrictive smoothness assumption.

Finally, we discuss the difference of the current paper and our prior works in \citep{doan2024fast,zeng2024fast}, which study the special setting $N=2$ under i.i.d. samples. First, this paper extends our prior works to the general setting with any $N$. Moving from $N=2$ to $N\geq 3$ creates a critical technical challenge -- with $N\geq 3$ mid-level variables $\{\theta_i\}_{i=2,\cdots,N-1}$ are coupled with faster and slower moving variables at the same time, making it harder for them to be decoupled by gradient averaging. We overcome the challenge through an improved bound on $x_i^{[t]}$ (which measures the convergence of $\theta_i^{[t]}$ to a proper learning target introduced later) by a function of $\{x_j^{[t]}\}_{j<i}$ and $\{x_j^{[t]}\}_{j>i}$. We also perform a mathematical investigation of the root cause of the complexity improvement. An explanation of why \texttt{A-MSA} enjoys a faster convergence rate is missing in our prior papers. This work fills in the gap by contrasting an important intermediate result obtained under \texttt{A-MSA} and \texttt{MSA} and pinpointing the terms that get improved due to the averaged operator estimation step. Finally, this paper  considers a more general sampling oracle -- we study non i.i.d samples generated from a time-varying Markov chain, a setting realistic for modeling applications in sequential decision making. To avoid complexity degradation due to Markovian samples, analyzing \texttt{A-MSA} requires carefully bounding the distance between the distribution of $X^{[t]}$ and a (time-varying) stationary distribution determined by $\vtheta^{[t]}$ by $\Ocal(\log t)$, whereas this distance is zero in expectation in the i.i.d. setting.

\section{Accelerated Multi-Time-Scale Stochastic Approximation Method}
The \texttt{MSA} method \eqref{alg:update_standardSA} updates each $\theta_i$ in the direction of the corresponding operator $F_i$ evaluated at the latest iterates. Our proposed method, formally presented in Algorithm \ref{alg:main}, modifies the update by introducing the auxiliary variables $f_i^{[t]}$ to estimate $\widebar{F}_i^{[t]}(\vtheta^{[t]})$ by forming a weighted average of the operator samples controlled by a parameter $\lambda^{[t]}\in(0,1)$. These auxiliary variables will then be used to update the decision variables. 
If we set $\lambda^{[t]}=1$ for all $t$, Algorithm~\ref{alg:main} reduces to \texttt{MSA} in \eqref{alg:update_standardSA}. However, a large $\lambda^{[t]}$ leads to a high-variance estimate, whereas if $\lambda^{[t]}$ is too small $f_i^{[t]}$ cannot track the moving target $\widebar{F}_i^{[t]}(\vtheta^{[t]})$. Decaying $\lambda^{[t]}$ appropriately with respect to $\alpha_i^{[t]}$, $f_{i}^t$ helps to control the variance of our estimates, which is crucial to derive the optimal convergence rate.

\begin{algorithm}
\caption{Accelerated Multi-Time-Scale Stochastic Approximation Algorithm}
\label{alg:main}
\begin{algorithmic}[1]
\STATE{\textbf{Initialization:} decision variables $\vtheta^{[0]}\in\Omega$, operator estimation variables $\{f_i^{[0]}\in\mathbb{R}^d\}_i$, initial sample $X^{[0]}$, step size sequences $\{\alpha_i^{[t]}\}_{i,t},\{\lambda^{[t]}\}_t$}
\FOR{$t=0,1,\cdots,T-1$}
\STATE{Decision variable update:
        \vspace{-5pt}
        \begin{align}
        \begin{aligned}
        \theta_i^{[t+1]} = \theta_i^{[t]} - \alpha_i^{[t]} f_i^{[t]},\;\forall i=1,\cdots,N,\qquad\quad \vtheta^{[t]}=\{\theta_1^{[t]},\cdots,\theta_N^{[t]}\}.
        \end{aligned}
        \label{alg:update_decision}
        \end{align}
        \vspace{-15pt}
        }
\STATE{Gradient estimation variable update:
\vspace{-5pt}
\begin{align}
\begin{aligned}
f_i^{[t+1]} = (1-\lambda^{[t]})f_i^{[t]} + \lambda^{[t]} F_i(\theta_1^{[t]},\cdots,\theta_N^{[t]},X^{[t]}).
\end{aligned}
\label{eq:update_auxiliary}
\end{align}}
\vspace{-25pt}
\STATE{Sample draw
\begin{align}
X^{[t+1]}\sim P_{\vtheta^{[t]}}(X^{[t]}).\label{alg:update_sample}
\end{align}}
\ENDFOR
\end{algorithmic}
\end{algorithm}

Note that in our algorithm the samples are drawn from the transition kernel $P_{\vtheta}$ parameterized by $\vtheta$ (a.k.a \eqref{alg:update_sample}). Thus, as $\vtheta^{[t]}$ changes over time, the samples are temporally dependent and Markovian with time-varying stationary distributions. 
To ensure the stable behavior of the Markovian samples we make the following assumptions. Our first assumption is on the fast mixing of the Markov chain generated under $P_{\vtheta}$. 
Given two distributions $u_1,u_2\in\Delta_{\Xcal}$, their total variation (TV) distance and the mixing time of the Markov chain are defined as
\begin{equation}
    \label{eq:TV_def}
    d_{\text{TV}}(u_1,u_2)=\frac{1}{2} \sup _{\nu: \Xcal \rightarrow[-1,1]}\left|\int \nu d u_1-\int \nu d u_2\right|.
\end{equation}

\begin{definition}
    \label{def:mixing_time}
    Given $\vtheta\in\Omega$, consider a Markov chain $\{X^{[t]}\}$ where $X^{[t]}\sim P_{\vtheta}(X^{[t-1]})$ with $\mu_{\vtheta}$ as the stationary distribution.  For any $a>0$, the mixing time of $\{X^{[t]}\}$ at level $a$ is defined as
    \[
        \tau_{\vtheta}(\alpha) = \min\{t\in\mathbb{N}:\sup_{X\in\Xcal}d_{\text{TV}}(\mathbb{P}(X^{[t]}=\cdot\mid X^{[0]}=X),\mu_{\vtheta})\leq a\}.
    \]
\end{definition}

\begin{assump}[Uniform Geometric Ergodicity]\label{assump:markov-chain}
    For any $\vtheta\in\Omega$, the Markov chain $\{X^{[t]}\}$ generated under $P_{\vtheta}$ has a unique stationary distribution $\mu_{\vtheta}$. In addition, there exist constants $m>0$ and $\rho\in (0,1)$ independent of $\vtheta$ such that
    \begin{equation*}
        \sup_{X\in\Xcal}d_{\text{TV}}(\mathbb{P}(X^{[t]}=\cdot\,|X^{[0]}=X),\mu_{\vtheta})\leq m\rho^t \text{ for all } t\geq 0. 
    \end{equation*}
\end{assump}
The geometric ergodicity assumption states that the distribution of a sample from the Markov chain approaches the stationary distribution geometrically fast, and is again standard in the literature on RL and SA under Markovian noise \citep{srikant2019finite,wu2020finite,zeng2024two}.

Denoting $\tau(a)\triangleq\sup_{\vtheta\in\mathbb{R}^d}\tau_{\vtheta}(a)$, this assumption implies that there exists a positive constant $C$ depending only on $\rho$ and $m$ such that 
\begin{equation}
    \tau(a) \leq C\log\left(1/a\right).\label{assump:mixing:tau}
\end{equation}
For convenience, we denote by $\tau_t\triangleq\tau(\lambda^{[t]})$ the mixing time corresponding to $\lambda^{[t]}$ in Algorithm \ref{alg:main}.

\begin{assump}
    \label{assump:tv_bound}
    Given two distributions $d,\hat{d}$ over $\Xcal$ and $\vtheta,\hat{\vtheta}\in\Omega$, let
    $X\sim d, X'\sim P_{\vtheta}( X)$ and $\hat{X}\sim\hat{d},\hat{X}'\sim P_{\hat{\vtheta}}(\hat{X})$.
    There exists a constant $L>0$ such that
    \begin{align}
        &d_{\text{TV}}(P(X'=\cdot), P(\hat{X}'=\cdot)) \leq d_{\text{TV}}(d, \hat{d}) + L\|\vtheta-\hat{\vtheta}\|.
        \label{assump:tv_bound:eq1}
    \end{align}
    In addition, the stationary distribution is Lipschitz in $\vtheta$
    \begin{align}
        d_{\text{TV}}(\mu_{\vtheta}, \mu_{\hat{\vtheta}})\leq L\|\vtheta-\hat{\vtheta}\|.
        \label{assump:tv_bound:eq2}
    \end{align}
\end{assump}
Assumption~\ref{assump:tv_bound} represents a regularity condition on $P_{\vtheta}$, which can be interpreted as a type of the triangle inequality for the TV distance defined over the transition kernels w.r.t different parameters. 
This assumption can be shown to hold in standard MDPs (See \citet[Lemma A1]{wu2020finite}).

We denote by $\Fcal_t=\{X^{[0]},X^{[1]},\cdots,X^{[t]}\}$ the filtration containing all randomness generated by Algorithm \ref{alg:main} up to time $t$.
\section{Main Results}

In this section, we establish the finite-sample complexity of \texttt{A-MSA} to the solution of \eqref{eq:obj} when the operators $F_{i}$ have a nested structure. For convenience, we refer to the first and last equations in \eqref{eq:obj} as the highest- and lowest-level problems, respectively, i.e., the problem levels increase from $N$ to 1. Levels here indicate the order in which the problem can be solved. 
With some abuse of notation, we denote $\vtheta_{1:N-1} = \{\theta_1,\cdots,\theta_{N-1}\}$. Given $\vtheta_{1:N-1}$, we use $y_N(\vtheta_{1:N-1})$ to denote the solution of the lowest level problem
\begin{align*}
\widebar{F}_N\Big(\vtheta_{1:N-1},y_N(\vtheta_{1:N-1})\Big)=0.
\end{align*}
Similarly, we denote by $\{y_{N-1}(\vtheta_{1:N-2}), y_N(\vtheta_{1:N-2})\}$ the solution of the last two levels corresponding to $\widebar{F}_{N-1},\widebar{F}_N$ given $\vtheta_{1:N-2}$
\begin{align}
\begin{aligned}
&\widebar{F}_{N-1}\big(\vtheta_{1:N-2},y_{N-1}(\vtheta_{1:N-2}), y_N(\vtheta_{1:N-2})\big)=0.\\
&\widebar{F}_{N}\big(\vtheta_{1:N-2},y_{N-1}(\vtheta_{1:N-2}), y_N(\vtheta_{1:N-2})\big)=0.
\end{aligned}
\label{eq:two_equations}
\end{align}
Conceptually, $y_{N-1}(\vtheta_{1:N-2})$ and $y_N(\vtheta_{1:N-2})$ represent the optimal solutions of the last and second last equation with respect to the last and second last decision variables when $\theta_1,\cdots,\theta_{N-2}$ are given.

Generalizing this line of discussion, we introduce $y_j(\vtheta_{1:i-1})$ to represent the optimal solution of the $j^{\text{th}}$ equation with respect to the $j^{\text{th}}$ variable when $\theta_{1},\cdots,\theta_{i-1}$ are given, for any $i\in[2,N]$ and $j\in[i,N]$. 
For $i=1$ (i.e. no $\theta$ is given), we additionally define $y_j(\emptyset)=\theta_j^{\star}\in\mathbb{R}^{d_k}$, which is the $j^{\text{th}}$ component of the optimal solution to \eqref{eq:obj}.
We introduce the aggregate notation
\begin{align*}
\ybf_{i:j}(\vtheta_{1:i-1}) \triangleq \{y_i(\vtheta_{1:i-1}),\cdots,y_j(\vtheta_{1:i-1})\}
\end{align*}
and note that $\ybf_{i:N}$ is the solution to the system
\begin{align*}
\widebar{F}_{j}\big(\vtheta_{1:i-1}, \ybf_{i:N}(\vtheta_{1:i-1})\big)\hspace{-1pt}=\hspace{-1pt}\widebar{F}_j\big(\vtheta_{1:i-1}, y_i(\vtheta_{1:i-1}),y_{i+1}(\vtheta_{1:i-1}),\cdots,y_N(\vtheta_{1:i-1})\big)\hspace{-1pt}=\hspace{-1pt}0,\,\forall j\in[i,N].
\end{align*}
%
We now introduce the main assumption that drives our analysis, referred to as the nested strong monotonicity of the operators $\widebar{F}_i$ corresponding to the increasing levels defined above.
\begin{assump}[Strong Monotonicity]\label{assump:stronglymonotone}
There exists a constant $\delta>0$ such that for all $i=1,\cdots,N$ and $\theta_1\in\mathbb{R}^{d_1},\cdots,\theta_{i-1}\in\mathbb{R}^{d_{i-1}}$, $\theta_i,\theta_i' \in \mathbb{R}^{d_i}$
\begin{align*}
    &\Big\langle \widebar{F}_i(\vtheta_{1:i-1},\theta_i,\ybf_{i+1:N}(\vtheta_{1:i-1},\theta_{i}) - \widebar{F}_i(\vtheta_{1:i-1},\theta_i',\ybf_{i+1:N}(\vtheta_{1:i-1},\theta_{i}'), \theta_i-\theta_i'\Big\rangle \geq\delta \|\theta_i-\theta_i'\|^2.
\end{align*}
\end{assump}
Assumption~\ref{assump:stronglymonotone} represents a nested variant of strong monotonicity for the operators $\widebar{F}_{i}$, which does not require each $\widebar{F}_{i}$ being strongly monotone w.r.t to all of its variables. It states that given $\vtheta_{i-1}$, $\widebar{F}_i$ is strongly monotone w.r.t $\theta_i$ when all lower level decision variables are at the corresponding optimal solutions. When $N=2$, the assumption reduces to the same condition made in the existing literature \citep{doan2022nonlinear,shen2022single,han2024finite} and can be verified to hold in RL applications including the gradient-based TD learning \citep{xu2019two} and distributed TD learning learning\footnote{To model the distributed TD learning in the two-time-scale SA framework, the outer loop equation enforces the consensus among agents whereas in the inner loop each agent learns its local value function.}. As an implication of the assumption, $y_j(\vtheta_{1:i-1})$ is always unique for any $i,j$. Without any loss of generality, we assume $\delta\leq 1$, which makes it convenient for us to simplify terms.

\begin{assump}[Lipschitz Continuity and Boundedness]\label{assump:Lipschitz}
There exists positive finite constants $L$ and $B$ such that for all $\vtheta,\vtheta'\in\Omega$
\begin{align}
    \|F_i(\vtheta,X)-F_i(\vtheta',X)\|&\leq L\sum_{j=1}^{N}\|\theta_j-\theta_j'\|,\quad\forall i,X,\label{assump:Lipschitz:eq1}\\
    \|y_j(\vtheta_{1:i-1})-y_j(\vtheta_{1:i-1}')\|&\leq L\sum_{i'=1}^{i-1}\|\theta_{i'}-\theta_{i'}'\|,\quad\forall i,j,\label{assump:Lipschitz:eq2}\\
    \|y_i(\vtheta_{1:i-1})\|&\leq B,\quad\forall i.\label{assump:Lipschitz:eq3}
\end{align}
\end{assump}
Eqs.~\eqref{assump:Lipschitz:eq1} and \eqref{assump:Lipschitz:eq2} are standard Lipschitz continuity assumptions on the operator $F_i$ and learning target \citep{doan2022nonlinear,zeng2024two,shen2022single}. Eq.~\ref{assump:Lipschitz:eq3} assumes that the learning targets are bounded and guarantees algorithm stability in the presence time-varying Markovian noise. This assumption is made in \citet{kaledin2020finite,zeng2024two} when $N = 2$.
Without any loss of generality, we assume $L\geq\max_{i\in\{1,\cdots,N\},X\in\Xcal}\|F_i(0,X)\|$, which is finite due to the compactness of $\Xcal$. Assumption~\ref{assump:Lipschitz} implies that the operator $F_i$ can be bounded affinely by the decision variable
\begin{align}\label{eq:bounded_F}
\|F_i(\vtheta,X)\|\leq L(\sum_{j=1}^{N}\|\theta_j\|+1),\quad\forall i\in\{1,\cdots,N\},\vtheta\in\Omega.
\end{align}
We use the same constant $L$ in Assumptions~\ref{assump:tv_bound} and \ref{assump:Lipschitz} to simplify the notation.

\noindent\textbf{Convergence Metric.} We quantity the algorithm convergence through the residuals
\begin{align}
\begin{gathered}
x_i^{[t]}=\theta_i^{[t]}-y_i(\vtheta_{1:i-1}^{[t]}),\qquad \Delta f_i^{[t]}=f_i^{[t]}-\widebar{F}_i(\vtheta^{[t]}),
\end{gathered}\label{eq:residuals}
\end{align}
where $x_i^{[t]}$ is the optimal residual, measuring the distance of decision variable $\theta_i^{[t]}$ to the learning target solution introduced above; whereas $\Delta f_i^{[t]}$ is the estimation residual, capturing the quality of the estimates of the operator $\widebar{F}_i$. Note that if $x_i^{[t]}=0$ for all $i$, then $\vtheta^{[t]} = \vtheta^{\star}$, the desired solution.

Finally, we will study the convergence of Algorithm \ref{alg:main} under the following choice of step sizes
\begin{align}
\lambda^{[t]} = \frac{c_{\lambda}}{t+h+1},\quad \alpha_i^{[t]} = \frac{c_i}{t+h+1},\label{eq:step_sizes}
\end{align}
where $\lambda^{[t]} \geq \alpha_i^{[t]}$ as the operator estimates are implemented at a ``faster" time scale than the updates of decision variables. In particular, these step sizes satisfies the following conditions for all $t\geq 0$
\begin{align}
    c_1&=\frac{32}{\delta},\quad \lambda^{[t]}\leq\frac{1}{4}, \quad \tau_t^2\lambda_{t-\tau_t}\leq\frac{1}{8DN^3}, \quad\alpha_i^{[t]}\leq\min\Big\{\frac{\delta^2}{80N^8 L^6},\frac{\delta}{40N^5 L^6},\frac{2}{5NL^2},\frac{1}{\delta}\Big\},\notag\\
    \frac{\alpha_i^{[t]}}{\lambda^{[t]}}&\leq\min\Big\{\frac{1}{8(DN^3+3/\delta+L)},\frac{\delta}{32DN^3},\frac{\delta}{32(\frac{9N^4 L^6}{\delta}+8N^3 L^3)},\frac{16}{\delta}\Big\},\notag\\
    \frac{\alpha_{i-1}^{[t]}}{\alpha_i^{[t]}}&\leq\frac{\delta}{16}\Big(\frac{9N^3 L^3}{2}+\frac{4N^6 L^6}{\delta}\Big)^{-1},\quad\frac{(\alpha_{i}^{[t]})^2}{\alpha_1^{[t]}}\leq\min\Big\{\frac{\delta^{3/2}}{64N^7},\frac{8}{5N^3 L^3}\Big\},
\label{eq:step_sizes_condition}
\end{align}
where $D$ is an absolute constant depending on $L$, $B$, $\rho$, and $m$.
We note that there always exist constants $c_{\lambda},c_{i},h$ such that the step sizes in \eqref{eq:step_sizes}  satisfy \eqref{eq:step_sizes_condition}. 

For our analysis, we consider the following Lyapunov function
\begin{align}
V^{[t]}=\sum_{i=1}^{N}\mathbb{E}[\|x_i^{[t]}\|^2+\|\Delta f_i^{[t]}\|^2].\label{thm:main:Lyapunov}
\end{align}
Under the choice of step sizes above, the main result of this paper is presented in the following theorem, where we characterize the convergence complexity of Algorithm \ref{alg:main}

\begin{thm}\label{thm:main}
Under Assumptions~\ref{assump:markov-chain}--\ref{assump:Lipschitz}, the iterates generated by Algorithm~\ref{alg:main} satisfy for all $t\geq\tau_t$
\begin{align*}
V^{[t+1]}\leq\frac{h^2 V^{[\tau_t]}}{(t+h+1)^2}+\widetilde{\Ocal}\left(\frac{1}{t+1}\right).
\end{align*}
\end{thm}
Theorem~\ref{thm:main} shows that the last iterate of the \texttt{A-MSA} algorithm converges to the unique solution of \eqref{eq:obj} in the mean-squared sense at an optimal rate $\widetilde{\Ocal}(1/t)$. 
As the algorithm draws exactly one sample in each iteration, this rate translates to a sample complexity of $\widetilde{\Ocal}(\epsilon^{-1})$ for finding an $\epsilon$-optimal solution. The same complexity has been shown to be achievable by the \texttt{MSA} algorithm under an additional smoothness assumption on $F_i$ \citep{shen2022single}. We match the rate without requiring this restrictive assumption.
In the absence of the smoothness condition, our rate is order-optimal (up to an logarithm factor) and significantly improves over the best-known rate of \texttt{MSA}, which is $\widetilde{\Ocal}(1/t^{2/3})$ as established in \citet{hong2023two,zeng2024two,han2024finite} when $N=2$. In Section~\ref{sec:standardSA}, we further show that the convergence rate of \texttt{MSA} is $\widetilde{\Ocal}(1/t^{1/2})$ when $N=3$ and deteriorates as $N$ increases. This observation highlights the advantage of our proposed algorithm where its convergence complexity is optimal for arbitrary $N$ .

\subsection{Proof of Theorem~\ref{thm:main}}

We introduce the following technical lemmas, which set up important intermediate steps to derive the results in Theorem~\ref{thm:main}. These
lemmas establishes a useful upper bound of $\|\Delta f_i^{[t]}\|^2$ and $\|x_i^{[t]}\|^2$ in expectation, which shows the dependence of these two quantities. 
\looseness=-1

\begin{lem}\label{lem:bound_Deltaf}
Under Assumptions~\ref{assump:markov-chain}-\ref{assump:Lipschitz}, the iterates of Algorithm~\ref{alg:main} satisfy for all $t\geq\tau_t$
\begin{align*}
\mathbb{E}[\|\Delta f_i^{[t+1]}\|^2]
&\leq (1-\lambda^{[t]})\mathbb{E}[\|\Delta f_i^{[t]}\|^2]-\frac{\lambda^{[t]}}{4}\mathbb{E}[\|\Delta f_i^{[t]}\|^2]+DN^2\tau_t^2\lambda^{[t]}\lambda^{[t-\tau_t]}\notag\\
&\hspace{20pt}+D N^2 \tau_t^2  (\lambda^{[t]}\lambda^{[t-\tau_t]}+\frac{(\alpha_i^{[t]})^2}{\lambda^{[t]}})\big(\sum_{j=1}^{N}\|x_j^{[t]}\|^2+\sum_{j=1}^{N}\|\Delta f_j^{[t]}\|^2\big),
\end{align*}
where $D$ is a constant depending only (polynomially) on the problem constants.
\end{lem}

\begin{lem}\label{lem:bound_x}
Under Assumptions~\ref{assump:markov-chain}-\ref{assump:Lipschitz}, the iterates of Algorithm~\ref{alg:main} satisfy for all $t$
\begin{align}
\|x_i^{[t+1]}\|^2&\leq \|x_i^{[t]}\|^2-\frac{\delta\alpha_i^{[t]}}{4}\|x_i^{[t]}\|^2+\sum_{j=1}^{i-1}\frac{\delta\alpha_j^{[t]}}{8N}\|x_j^{[t]}\|^2\notag\\
&\hspace{20pt}+\left(\frac{9N^3 L^6}{\delta}+8N^2 L^3\right)\alpha_i^{[t]}\sum_{j=i+1}^{N}\|x_j^{[t]}\|^2+(\frac{3}{\delta}+L)\alpha_i^{[t]}\sum_{j=1}^{N}\|\Delta f_j^{[t]}\|^2.\label{lem:bound_x:eq1}
\end{align}
\end{lem}
In \eqref{lem:bound_x:eq1}, the terms $\|x_j^{[t]}\|^2$ for $j<i$ and $j>i$ are scaled by different step sizes ($\alpha_j^{[t]}$ versus $\alpha_i^{[t]}$). This subtle difference  is the key to establish the optimal complexity of Algorithm \ref{alg:main}. 

It is worth noting that if the iterates $\theta_i^{[t]}$ are generated by the \texttt{MSA} algorithm, a per-iteration analysis of $x_i^{[t]}$ similar to Lemma~\ref{lem:bound_x} can be established. We make the analysis and point out the distinctions that allow \texttt{A-MSA} to converge faster later in Section~\ref{sec:standardSA} (Lemma~\ref{lem:conv_x_standardSA}).

We defer the proofs of the lemmas to the appendix but point out that to prove Lemma~\ref{lem:bound_x} the following bound on $\|\theta_i^{[t+1]}-\theta_i^{[t]}\|$ provides an important tool for us to control the stability of the iterates. The implication of Lemma~\ref{lem:Lipschitz} is that the distance $\|\theta_i^{[t+1]}-\theta_i^{[t]}\|$ converges to zero faster than $\Ocal(\alpha_i^{[t]})$ if $\|\Delta f_i^{[t]}\|$ and $\sum_{j=i}^{N}\|x_j^{[t]}\|$ both converge.
\begin{lem}\label{lem:Lipschitz}
Under Assumptions~\ref{assump:markov-chain}-\ref{assump:Lipschitz}, the iterates of Algorithm~\ref{alg:main} satisfy for all $i,t$
\begin{align*}
\|\theta_i^{[t+1]}-\theta_i^{[t]}\|\leq\alpha_i^{[t]}\big(\|\Delta f_i^{[t]}\|+NL^2\sum_{j=i}^{N}\|x_j^{[t]}\|\big).
\end{align*}
\end{lem}

\noindent\textbf{Proof (of Theorem~\ref{thm:main}).}

The proof follows in a straightforward manner by combining the bounds in Lemmas~\ref{lem:bound_Deltaf} and \ref{lem:bound_x} and choosing the correct step sizes. Recall the Lyapunov function \eqref{thm:main:Lyapunov}
\begin{align}
V^{[t+1]}&=\sum_{i=1}^{N}\mathbb{E}[\|\Delta f_i^{[t+1]}\|^2+\|x_i^{[t+1]}\|^2]\notag\\
&\leq (1-\lambda^{[t]})\sum_{i=1}^{N}\mathbb{E}[\|\Delta f_i^{[t]}\|^2]-\frac{\lambda^{[t]}}{4}\sum_{i=1}^{N}\mathbb{E}[\|\Delta f_i^{[t]}\|^2]+DN^2\tau_t^2\lambda^{[t]}\lambda^{[t-\tau_t]}\notag\\
&\hspace{20pt}+D N^2 \tau_t^2  \sum_{i=1}^N(\lambda^{[t]}\lambda^{[t-\tau_t]}+\frac{(\alpha_i^{[t]})^2}{\lambda^{[t]}})\big(\sum_{j=1}^{N}\|x_j^{[t]}\|^2+\sum_{j=1}^{N}\|\Delta f_j^{[t]}\|^2\big)\notag\\
&\hspace{20pt}+\sum_{i=1}^{N}\mathbb{E}[\|x_i^{[t]}\|^2]-\sum_{i=1}^{N}\frac{\delta\alpha_i^{[t]}}{4}\mathbb{E}[\|x_i^{[t]}\|^2] +\sum_{i=1}^{N}\sum_{j=1}^{i-1}\frac{\delta\alpha_j^{[t]}}{8N}\mathbb{E}[\|x_j^{[t]}\|^2]\notag\\
&\hspace{20pt}+\sum_{i=1}^{N}\left(\frac{9N^3 L^6}{\delta}+8N^2 L^3\right)\alpha_i^{[t]}\sum_{j=i+1}^{N}\mathbb{E}[\|x_j^{[t]}\|^2]+(\frac{3}{\delta}+L)\sum_{i=1}^{N}\alpha_i^{[t]}\sum_{j=1}^{N}\mathbb{E}[\|\Delta f_j^{[t]}\|^2]\notag\\
&\leq (1-\lambda^{[t]})\sum_{i=1}^{N}\mathbb{E}[\|\Delta f_i^{[t]}\|^2]-\frac{\lambda^{[t]}}{4}\sum_{i=1}^{N}\mathbb{E}[\|\Delta f_i^{[t]}\|^2]+\widetilde{\Ocal}\Big(\frac{1}{(t+1)^2}\Big)\notag\\
&\hspace{20pt}+DN^3\alpha_N^{[t]}\sum_{i=1}^{N}\mathbb{E}[\|\Delta f_i^{[t]}\|^2] + DN^3 \sum_{i=1}^{N}\frac{(\alpha_i^{[t]})^2}{\lambda^{[t]}}\mathbb{E}[\|x_i^{[t]}\|^2]\notag\\
&\hspace{20pt}+DN^3\tau_t^2\lambda^{[t]}\lambda^{[t-\tau_t]}\sum_{i=1}^{N}\mathbb{E}[\|\Delta f_i^{[t]}\|^2] + DN^3 \tau_t^2\lambda^{[t]}\lambda^{[t-\tau_t]}\sum_{i=1}^{N}\mathbb{E}[\|x_i^{[t]}\|^2]\notag\\
&\hspace{20pt}+\sum_{i=1}^{N}\mathbb{E}[\|x_i^{[t]}\|^2]-\sum_{i=1}^{N}\frac{\delta\alpha_i^{[t]}}{4}\mathbb{E}[\|x_i^{[t]}\|^2] +\sum_{i=1}^{N}\frac{\delta\alpha_i^{[t]}}{8}\mathbb{E}[\|x_i^{[t]}\|^2]\notag\\
&\hspace{20pt}+\left(\frac{9N^4 L^6}{\delta}+8N^3 L^3\right)\sum_{i=2}^{N}\alpha_{i-1}^{[t]}\mathbb{E}[\|x_i^{[t]}\|^2]+(\frac{3}{\delta}+L)\alpha_N^{[t]}\sum_{i=1}^{N}\mathbb{E}[\|\Delta f_i^{[t]}\|^2]\notag\\
&\leq (1-\lambda^{[t]})\sum_{i=1}^{N}\mathbb{E}[\|\Delta f_i^{[t]}\|^2]+(1-\frac{\delta\alpha_1^{[t]}}{16})\sum_{i=1}^{N}\mathbb{E}[\|x_i^{[t]}\|^2]+\widetilde{\Ocal}\Big(\frac{1}{(t+1)^2}\Big)\notag\\
&\hspace{20pt}-\sum_{i=1}^{N}\left(\frac{\lambda^{[t]}}{8}-(DN^3+\frac{3}{\delta}+L)\alpha_N^{[t]}\right)\mathbb{E}[\|\Delta f_i^{[t]}\|^2]\notag\\
&\hspace{20pt}-\sum_{i=1}^{N}\left(\frac{\delta\alpha_i^{[t]}}{16}-\frac{DN^3 (\alpha_i^{[t]})^2}{\lambda^{[t]}}-\left(\frac{9N^4 L^6}{\delta}+8N^3 L^3\right)\alpha_{i-1}^{[t]}\right)\mathbb{E}[\|x_i^{[t]}\|^2],\label{alg:main:proof_eq1}
\end{align}
where the first inequality plugs in the bounds from Lemmas~\ref{lem:bound_Deltaf} and \ref{lem:bound_x}, and in the last inequality we have artificially introduced $\alpha_0^{[t]}=0$ and used the relation $DN^3\tau_t^2\lambda^{[t-\tau_t]}\leq1/8$. Note that the last two terms of \eqref{alg:main:proof_eq1} are non-positive under the step size condition $\alpha_i^{[t]}\leq\frac{\lambda^{[t]}}{8}(DN^3+\frac{3}{\delta}+L)^{-1}$, $\alpha_i^{[t]}\leq\frac{\delta\lambda^{[t]}}{32DN^3}$, and $\alpha_{[i-1]}^{[t]}\leq\frac{\delta\alpha_i^{[t]}}{32}(\frac{9N^4 L^6}{\delta}+8N^3 L^3)^{-1}$. As a result, we have
\begin{align*}
V^{[t+1]}&\leq\max\{1-\lambda^{[t]},1-\frac{\delta\alpha_1^{[t]}}{16}\}V^{[t]}+\widetilde{\Ocal}\Big(\frac{1}{(t+1)^2}\Big)\notag\\
&\leq(1-\frac{\delta\alpha_1^{[t]}}{16})V^{[t]}+\widetilde{\Ocal}\Big(\frac{1}{(t+1)^2}\Big)\notag\\
&\leq (1-\frac{2}{t+h+1})V^{[t]}+\widetilde{\Ocal}\Big(\frac{1}{(t+1)^2}\Big),
\end{align*}
due to $\alpha_{1}^{[t]}\leq\frac{16\lambda^{[t]}}{\delta}$ and $c_1=\frac{32}{\delta}$. Multiplying by $(t+h+1)^2$ on both sides, we have
\begin{align*}
(t+h+1)^2 V^{[t+1]}&\leq(t+h+1)(t+h-1)V^{[t]}+\widetilde{\Ocal}(1)\notag\\
&\leq (t+h)^2 V^{[t]}+\widetilde{\Ocal}(1)\notag\\
&\leq h^2 V^{[\tau_t]}+\widetilde{\Ocal}(t+1),
\end{align*}
which when dividing by $(t+h+1)^2$ gives the desired result. 


\qed
\subsection{Standard Multi-Time-Scale Stochastic Approximation}\label{sec:standardSA}

In this section, we present the finite-time complexity of the \texttt{MSA} algorithm derived using the techniques in \citet{hong2023two,zeng2024two}, to highlight its gap from the optimal rate $\widetilde{\Ocal}(1/t)$ established in Theorem~\ref{thm:main}. We start with the case of $N=3$, where we want to solve
\begin{align*}
     \widebar{F}_{i}(\vtheta)=0,\quad \forall i = 1,2,3.
\end{align*}
For simplicity, we assume (just in this section) that the stationary distribution $\mu_{\vtheta}=\mu$ is independent of $\vtheta$ and that each $X^{[t]}$ is an i.i.d. sample from $\mu$. We additionally assume that $F_i$ has uniformly bounded energy as follows. Note that this is not a necessary assumption but considered just to simplify the analysis of \texttt{MSA}. Also note that we do not make this assumption in the analysis of the proposed \texttt{A-MSA} algorithm.

\begin{assump}\label{assump:bounded}
There exists a constant $D>0$ such that for all $i=1,2,3$, $\theta_i \in \mathbb{R}^{d_i}$, and $X\in\Xcal$
\begin{align*}
\| F_i(\theta_1,\theta_2,\theta_3,X)\|\leq D.
\end{align*}
\end{assump}

Recall the update rule of \texttt{MSA} in \eqref{alg:update_standardSA}.  We measure the convergence by $x_i^{[t]}$, the same as the metric \eqref{eq:residuals} used for analyzing \texttt{A-MSA}.
\begin{align*}
x_i^{[t]}=\theta_i^{[t]}-y_i(\vtheta_{1:i-1}^{[t]}).
\end{align*}
We next establish the per-iteration bound for $x_i^{t}$, where the analysis can be found in Appendix~\ref{sec:proof_lem:conv_x_standardSA}. 
\begin{lem}\label{lem:conv_x_standardSA}
We artificially define $\alpha_0^{[t]}=0$ for all $t$.
Under Assumptions~\ref{assump:stronglymonotone}-\ref{assump:bounded}, the iterates of the \texttt{MSA} algorithm in \eqref{alg:update_standardSA} satisfy for all $t\geq0$
\begin{align*}
\mathbb{E}[\|x_i^{[t+1]}\|^2]
&\leq \mathbb{E}[\|x_i^{[t]}\|^2]-\frac{\delta\alpha_i^{[t]}}{4}\mathbb{E}[\|x_i^{[t]}\|^2]+\frac{8N^3 L^6\alpha_i^{[t]}}{\delta}\sum_{j=i+1}^{N}\mathbb{E}[\|x_j^{[t]}\|^2]\\
&\hspace{120pt}+8N^2 L^2 D^2 (\alpha_{i}^{[t]})^2+\frac{N^2 L^2 D^2 (\alpha_{i-1}^{[t]})^2}{\delta\alpha_i^{[t]}}.
\end{align*}
\end{lem}
It is worth noting that Lemma~\ref{lem:conv_x_standardSA} holds for any $N$. We apply this lemma to the case $N=3$ and consider the Lyapunov function 
\begin{align}
V^{[t]}=\mathbb{E}[\|x_1^{[t]}\|^2+v_2^{[t]}\|x_2^{[t]}\|^2+v_3^{[t]}\|x_3^{[t]}\|^2],\label{eq:standardSA_lyapunov}
\end{align}
where $v_2^{[t]}=\frac{1728 L^6\alpha_1^{[t]}}{\delta^2\alpha_2^{[t]}}$ and $v_3^{[t]}=\frac{8\alpha_1^{[t]}}{\delta\alpha_3^{[t]}}(\frac{216 L^6}{\delta}+\frac{373248L^{12}}{\delta^3})$.

\begin{thm}\label{thm:standardSA}
Let the step sizes be
\begin{align}
\alpha_1^{[t]}=\frac{c_1}{t+h+1},\quad \alpha_2^{[t]}=\frac{c_2}{(t+h+1)^{3/4}},\quad \alpha_2^{[t]}=\frac{c_3}{(t+h+1)^{1/2}},\label{thm:standardSA:step_size}
\end{align}
with proper choices of the constants $c_1,c_2,c_3$, and $h$.
Under Assumptions~\ref{assump:stronglymonotone}-\ref{assump:bounded}, the iterates of the \texttt{MSA} algorithm in \eqref{alg:update_standardSA} with i.i.d. samples $X^{[t]}\sim\mu$ satisfy for all $t\geq0$
\begin{align*}
V^{[t+1]}\leq\Ocal\left(\frac{V^{[0]}}{(t+\tau+1)^{1/2}}\right).
\end{align*}
\end{thm}

\textbf{Proof.}

With $N=3$, Lemma~\ref{lem:conv_x_standardSA} implies
\begin{align*}
\mathbb{E}[\|x_1^{[t+1]}\|^2]&\hspace{-2pt}\leq\hspace{-2pt} (1\hspace{-2pt}-\hspace{-2pt}\frac{\delta\alpha_1^{[t]}}{4})\mathbb{E}[\|x_1^{[t]}\|^2]+\frac{216 L^6\alpha_1^{[t]}}{\delta}\mathbb{E}[\|x_2^{[t]}\|^2+\|x_3^{[t]}\|^2]+72 L^2 D^2 (\alpha_1^{[t]})^2,\\
\mathbb{E}[\|x_2^{[t+1]}\|^2]&\hspace{-2pt}\leq\hspace{-2pt} (1\hspace{-2pt}-\hspace{-2pt}\frac{\delta\alpha_2^{[t]}}{4})\mathbb{E}[\|x_2^{[t]}\|^2]\hspace{-2pt}+\hspace{-2pt}\frac{216 L^6\alpha_2^{[t]}}{\delta}\mathbb{E}[\|x_3^{[t]}\|^2]\hspace{-2pt}+\hspace{-2pt}72 L^2 D^2 (\alpha_{2}^{[t]})^2\hspace{-2pt}+\hspace{-2pt}\frac{9 L^2 D^2 (\alpha_1^{[t]})^2}{\delta\alpha_2^{[t]}},\\
\mathbb{E}[\|x_3^{[t+1]}\|^2]&\hspace{-2pt}\leq\hspace{-2pt} (1\hspace{-2pt}-\hspace{-2pt}\frac{\delta\alpha_3^{[t]}}{4})\mathbb{E}[\|x_3^{[t]}\|^2]+72 L^2 D^2 (\alpha_3^{[t]})^2+\frac{9 L^2 D^2 (\alpha_2^{[t]})^2}{\delta\alpha_3^{[t]}}.
\end{align*}

Recall the Lyapunov function 
$V^{[t]}=\mathbb{E}[\|x_1^{[t]}\|^2+v_2^{[t]}\|x_2^{[t]}\|^2+v_3^{[t]}\|x_3^{[t]}\|^2]$. Note that $v_2^{[t+1]}\leq v_2^{[t]}$ and $v_3^{[t+1]}\leq v_3^{[t]}$ since $\alpha_1^{[t]}$ decays faster than $\alpha_2^{[t]},\alpha_3^{[t]}$. We have
\begin{align}
V^{[t+1]}&=\mathbb{E}[\|x_1^{[t+1]}\|^2+v_2^{[t+1]}\|x_2^{[t+1]}\|^2+v_3^{[t+1]}\|x_3^{[t+1]}\|^2]\notag\\
&\leq \mathbb{E}[\|x_1^{[t+1]}\|^2+v_2^{[t]}\|x_2^{[t+1]}\|^2+v_3^{[t]}\|x_3^{[t+1]}\|^2]\notag\\
&\leq(1-\frac{\delta\alpha_1^{[t]}}{4})\mathbb{E}[\|x_1^{[t]}\|^2]+\frac{216 L^6\alpha_1^{[t]}}{\delta}\mathbb{E}[\|x_2^{[t]}\|^2+\|x_3^{[t]}\|^2]+72 L^2 D^2 (\alpha_{1}^{[t]})^2\notag\\
&\hspace{20pt}+v_2^{[t]}(1-\frac{\delta\alpha_2^{[t]}}{4})\mathbb{E}[\|x_2^{[t]}\|^2]+\frac{216 L^6\alpha_2^{[t]}v_2^{[t]}}{\delta}\mathbb{E}[\|x_3^{[t]}\|^2]+72 L^2 D^2 (\alpha_{2}^{[t]})^2 v_2^{[t]}\notag\\
&\hspace{20pt}+\frac{9 L^2 D^2 (\alpha_1^{[t]})^2 v_2^{[t]}}{\delta\alpha_2^{[t]}}+v_3^{[t]}(1-\frac{\delta\alpha_3^{[t]}}{4})\mathbb{E}[\|x_3^{[t]}\|^2]\notag\\
&\hspace{20pt}+72 L^2 D^2 (\alpha_3^{[t]})^2 v_3^{[t]}+\frac{9 L^2 D^2 (\alpha_2^{[t]})^2 v_3^{[t]}}{\delta\alpha_3^{[t]}}\notag\\
&=(1-\frac{\delta\alpha_1^{[t]}}{4})\mathbb{E}[\|x_1^{[t]}\|^2+v_2^{[t]}\|x_2^{[t]}\|^2+v_3^{[t]}\|x_3^{[t]}\|^2]+72 L^2 D^2 (\alpha_{1}^{[t]})^2\notag\\
&\hspace{20pt}+\frac{124416 L^8 D^2 \alpha_{1}^{[t]}\alpha_{2}^{[t]}}{\delta^2}+\frac{15552 L^8 D^2 (\alpha_1^{[t]})^3}{\delta^3(\alpha_2^{[t]})^2}\notag\\
&\hspace{20pt}+\left(\frac{572L^2 D^2\alpha_1^{[t]}\alpha_3^{[t]}}{\delta}+\frac{9L^2 D^2 \alpha_1^{[t]}(\alpha_2^{[t]})^2}{\delta^2(\alpha_3^{[t]})^2}\right)(\frac{216 L^6}{\delta}+\frac{373248L^{12}}{\delta^3})\notag\\
&\hspace{20pt}+(\frac{\delta \alpha_1^{[t]} v_2^{[t]}}{4}-\frac{\delta \alpha_2^{[t]} v_2^{[t]}}{4}+\frac{216 L^6\alpha_1^{[t]}}{\delta})\mathbb{E}[\|x_2^{[t]}\|^2]\notag\\
&\hspace{20pt}+(\frac{\delta \alpha_1^{[t]} v_3^{[t]}}{4}-\frac{\delta \alpha_3^{[t]} v_3^{[t]}}{4}+\frac{216 L^6\alpha_1^{[t]}}{\delta}+\frac{216 L^6\alpha_2^{[t]}v_2^{[t]}}{\delta})\mathbb{E}[\|x_3^{[t]}\|^2]\notag\\
&\leq(1-\frac{\delta\alpha_1^{[t]}}{4})V^{[t]}+72 L^2 D^2 (\alpha_{1}^{[t]})^2+\frac{124416 L^8 D^2 \alpha_{1}^{[t]}\alpha_{2}^{[t]}}{\delta^2}+\frac{15552 L^8 D^2 (\alpha_1^{[t]})^3}{\delta^3(\alpha_2^{[t]})^2}\notag\\
&\hspace{20pt}+\left(\frac{572L^2 D^2\alpha_1^{[t]}\alpha_3^{[t]}}{\delta}+\frac{9L^2 D^2 \alpha_1^{[t]}(\alpha_2^{[t]})^2}{\delta^2(\alpha_3^{[t]})^2}\right)(\frac{216 L^6}{\delta}+\frac{373248L^{12}}{\delta^3}),\label{thm:standardSA:proof_eq1}
\end{align}
where the last inequality is a result of 
\begin{align*}
&\frac{\delta \alpha_1^{[t]} v_2^{[t]}}{4}-\frac{\delta \alpha_2^{[t]} v_2^{[t]}}{4}+\frac{216 L^6\alpha_1^{[t]}}{\delta}=\frac{\delta \alpha_1^{[t]} v_2^{[t]}}{4}-\frac{\delta \alpha_2^{[t]} v_2^{[t]}}{8}\leq0,\\
&\frac{\delta \alpha_1^{[t]} v_3^{[t]}}{4}-\frac{\delta \alpha_3^{[t]} v_3^{[t]}}{4}+\frac{216 L^6\alpha_1^{[t]}}{\delta}+\frac{216 L^6\alpha_2^{[t]}v_2^{[t]}}{\delta}\\
&=\frac{\delta \alpha_1^{[t]} v_3^{[t]}}{4}-\frac{\delta \alpha_3^{[t]} v_3^{[t]}}{4}+\frac{216 L^6\alpha_1^{[t]}}{\delta}+\frac{373248L^{12}\alpha_1^{[t]}}{\delta^3}=\frac{\delta \alpha_1^{[t]} v_3^{[t]}}{4}-\frac{\delta \alpha_3^{[t]} v_3^{[t]}}{8}\leq0,
\end{align*}
following if the step sizes satisfy $\alpha_1^{[t]}\leq\frac{\alpha_2^{[t]}}{2}\leq \frac{\alpha_3^{[t]}}{2}$.

Under the step size rule \eqref{thm:standardSA:step_size} with $c_1\geq\frac{8}{\delta}$, it can be easily derived that
\begin{align*}
V^{[t+1]}&\leq(1-\frac{2}{t+h+1})V^{[t]}+\Ocal(1/t^{3/2}),
\end{align*}
which gives
\begin{align*}
V^{[t+1]}\leq\Ocal(\frac{V^{[0]}}{(t+h+1)^{1/2}}).
\end{align*}
\qed

We now discuss how Theorem~\ref{thm:standardSA} can be extended to more than three equations and what the expected convergence rate should be. 
When $N=2$, \citet{zeng2024two} shows that the convergence rate is on the order of (up to a logarithm factor)
\begin{align}
\max\Big\{\alpha_1^{[t]},\alpha_2^{[t]},\Big(\frac{\alpha_1^{[t]}}{\alpha_2^{[t]}}\Big)^2\Big\}.\label{eq:convrate_stepsize_N2}
\end{align}
To optimize \eqref{eq:convrate_stepsize_N2} the step sizes should be selected as $\alpha_1^{[t]}\sim\frac{1}{t}, \alpha_2^{[t]}\sim\frac{1}{t^{2/3}}$, which leads to a convergence rate of $\widetilde{\Ocal}(t^{-2/3})$. When $N=3$, the convergence rate is on the order of
\begin{align*}
\max\Big\{\alpha_1^{[t]},\alpha_2^{[t]},\alpha_3^{[t]},\Big(\frac{\alpha_1^{[t]}}{\alpha_2^{[t]}}\Big)^2,\Big(\frac{\alpha_2^{[t]}}{\alpha_3^{[t]}}\Big)^2\Big\}.
\end{align*}
which dictates the choice of step sizes in \eqref{thm:standardSA:step_size}.

For general $N$, the convergence rate becomes on the order of
\begin{align}
\max\Big\{\alpha_1^{[t]},\cdots,\alpha_N^{[t]},\Big(\frac{\alpha_1^{[t]}}{\alpha_2^{[t]}}\Big)^2,\cdots,\Big(\frac{\alpha_{N-1}^{[t]}}{\alpha_N^{[t]}}\Big)^2\Big\}.\label{eq:convrate_stepsize_N}
\end{align}
It is straightforward to see that the optimal decay rate of \eqref{eq:convrate_stepsize_N} is $\widetilde{\Ocal}(t^{-2/(N+1)})$, which implies the sample complexity in Table~\ref{table:comparison}.

\subsection{Reflecting on Faster Convergence of \texttt{A-MSA}}\label{sec:comparison}

Lemmas~\ref{lem:bound_x} and \ref{lem:conv_x_standardSA} establish critical per-iteration convergence bounds on the errors in the decision variable for \texttt{A-MSA} and \texttt{MSA}, respectively. In this section, we contrast corresponding terms in the bounds to pinpoint the advantage of Lemma~\ref{lem:bound_x} over Lemma~\ref{lem:conv_x_standardSA}. In the analysis of $\|x_i^{[t+1]}\|^2$ under the two algorithms, the common/comparable terms include the following:
\begin{itemize}
\item The same term $\|y_i(\vtheta_{1:i-1}^{[t+1]})-y_i(\vtheta_{1:i-1}^{[t]})\|$ arises from the shift in learning target $y_i(\vtheta_{1:i-1}^{[t]})$ across iterations and appears in the proofs of both Lemmas \ref{lem:bound_x} and \ref{lem:conv_x_standardSA} (in Section~\ref{sec:proof_lem:bound_x} and \ref{sec:proof_lem:conv_x_standardSA}, respectively). Due to the Lipschitz continuity of $y_i$ imposed in Assumption~\ref{assump:Lipschitz}, controlling this term essentially requires handling
\begin{align}
\sum_{j=1}^{i-1}\|\theta_j^{[t+1]}-\theta_j^{[t]}\|.\label{eq:contrast_eq1}
\end{align}

Under \texttt{MSA}, there is currently no better bound on \eqref{eq:contrast_eq1} than $\sum_{j=1}^{i-1}\alpha_j^{[t]}\|F_i(\vtheta^{[t]},X^{[t]})\|$. As $\|F_i(\vtheta^{[t]},X^{[t]})\|$ is on the order of a constant even when $\vtheta^{[t]}$ approaches the optimal solution (due to the random sample $X^{[t]}$), the bound eventually becomes $\Ocal(\alpha_{i-1})$.

In comparison, under \texttt{A-MSA}, we can control \eqref{eq:contrast_eq1} with Lemma~\ref{lem:Lipschitz} and show 
\begin{align}
    \sum_{j=1}^{i-1}\|\theta_j^{[t+1]}-\theta_j^{[t]}\|\leq\Ocal\Big(\sum_{j=1}^{i-1}\alpha_j^{[t]}\|\Delta f_j^{[t]}\|+\sum_{j=1}^{i-1}\alpha_j^{[t]}\|x_j^{[t]}\|+\alpha_{i-1}^{[t]}\sum_{j=i}^{N}\|x_j^{[t]}\|\Big).\label{eq:contrast_eq2}
\end{align}
Since $\|\Delta f_j^{[t]}\|$ and $\|x_j^{[t]}\|$ themselves decay to zero as the iterations proceed, the convergence rate of \eqref{eq:contrast_eq2} is much faster than $\Ocal(\alpha_{i-1})$.

\noindent\textbf{Takeaway:} 
Our objective is to solve a highly coupled system \eqref{eq:obj}.
The discussion above highlights how the \texttt{A-MSA} algorithm effectively decouples the decision variables, in the sense that the inaccuracy of the upper-level variables has a reduced effect on the lower-level ones. Being able to decoupled the decision variable updates is crucial as it allows us to choose $\alpha_i^{[t]}$ across $i$ in a more independent manner.

\item $\|\Delta f_i^{[t]}\|^2=\|f_i^{[t]}-\widebar{F}_i(\vtheta^{[t]})\|^2$ is an error appearing in the proof of Lemma~\ref{lem:bound_x} for \texttt{A-MSA}, whereas $\|F_i(\vtheta^{[t]},X^{[t]})-\widebar{F}_i(\vtheta^{[t]})\|^2$ is the comparable term in the proof of Lemma~\ref{lem:conv_x_standardSA} for \texttt{MSA}. They capture the variance in the estimation of $\widebar{F}_i(\vtheta^{[t]})$. In the case of \texttt{MSA}, $\|F_i(\vtheta^{[t]},X^{[t]})-\widebar{F}_i(\vtheta^{[t]})\|^2$ does not decay to zero (again due to the randomness in $X^{[t]}$), which eventually becomes a bottleneck in achieving the optimal rate. In contrast, $\|\Delta f_i^{[t]}\|^2$ under \texttt{A-MSA} can be expected to converge to zero (under a careful analysis).

\noindent\textbf{Takeaway:} 
Compared to \texttt{MSA}, \texttt{A-MSA} provides a low-variance estimate of $\widebar{F}_i(\vtheta^{[t]})$. The variance decaying sufficiently fast is another key driver of the overall improved complexity.

\end{itemize}

\begin{remark}
At the end of the section, we point out an important additional advantage of the analysis of \texttt{A-MSA} over that of \texttt{MSA}. Observe that the Lyapunov function \eqref{eq:standardSA_lyapunov} considered in Theorem~\ref{thm:standardSA} weights the lower-level residuals $x_2^{[t]}$ and $x_3^{[t]}$ by $v_2^{[t]}$ and $v_3^{[t]}$. The weights have been carefully chosen to ensure the proper cancellation of errors across different levels. Choosing such weights, however, requires non-trivial efforts and can be inconvenient as $N$ goes up. 

In addition, we note that the weights $v_2^{[t]}$ and $v_3^{[t]}$ themselves are decaying sequences. Theorem~\ref{thm:standardSA} shows that $V^{[t]}$ decays with rate $\widetilde{\Ocal}(t^{-1/2})$, which only implies $\mathbb{E}[\|x_2^{[t]}\|^2]\leq\widetilde{\Ocal}(t^{-1/4})$ and $\mathbb{E}[\|x_3^{[t]}\|^2]\leq\widetilde{\Ocal}(1)$. In other words, the guaranteed convergence rate of lower-level variables does not match that of the highest level variable and may even be meaningless.
In contrast, the Lyapunov function \eqref{thm:main:Lyapunov} for the analysis of \texttt{A-MSA} simply combines $\|x_i^{[t]}\|^2$ over all levels without additional weights. As a result, the decision variables at all level are guaranteed to converge at $\widetilde{\Ocal}(t^{-1})$.
\end{remark}

\section{Motivating Applications}\label{sec:applications}

In this section we discuss how the framework \eqref{eq:obj} models the objective of a range of algorithms in RL and games. In some cases, the objective is structured such that Assumptions~\ref{assump:stronglymonotone}-\ref{assump:tv_bound} can be verified to hold, meaning that the theoretical guarantee immediately applies if the \texttt{A-MSA} algorithm is used. For the other problems, the assumptions may not hold, but we discuss how the analysis can be adapted to derive the complexity of \texttt{A-MSA} by using the specific problem structure.
All applications discussed are special cases of \eqref{eq:obj} with $N=3$, meaning that they cannot be covered by the bilevel framework in \citet{zeng2024two,han2024finite,doan2024fast,zeng2024fast}.\looseness=-1 

\subsection{Gradient-Based Temporal Difference learning with Momentum.}\label{sec:applications:TDC}

Gradient-based temporal difference (TD) learning algorithms including GTD2 and TDC are a popular class of methods that stably minimize a projected Bellman error under linear function approximation and off-policy samples \citep{sutton2008convergent,sutton2009fast}. Consider an infinite-horizon average-reward MDP characterized by $(\Scal,\Acal,\Pcal,r,\gamma)$. Here $\Scal$ and $\Acal$ are the action and state spaces. Each state $s$ is associated with a feature vector $\phi(s)\in\mathbb{R}^d$. $\Pcal:\Scal\times\Acal\rightarrow\Delta_{\Scal}$ is the transition probability kernel, with $\Pcal(s'\mid s,a)$ denoting the probability of the MDP transitioning from state $s$ to $s'$ under action $a$. $r:\Scal\times\Acal\rightarrow[0,1]$ is the reward function. $\gamma\in(0,1)$ is the discount factor. Given a fixed policy $\pi$ to evaluate, our objective is to learn a value function parameter $\theta\in\mathbb{R}^d$ such that $\phi(s)^{\top}\theta$ approximates a value function $V^{\pi}(s)$ for every state $s\in\Scal$. The TDC algorithm formulates the objective with a system of two equations on $\theta$ and an auxiliary variable $\omega\in\mathbb{R}^d$, i.e., $\forall s\in\Scal$ we aim to solve
\begin{gather}
\begin{aligned}
&\mathbb{E}_{a\sim\pi(\cdot\mid s),s'\sim\Pcal(\cdot\mid s,a))}[r(s,a)+\gamma\phi(s')^{\top}\theta\phi(s)-\phi(s)^{\top}\theta\phi(s)-\gamma\phi(s')\phi(s)^{\top}\omega]=0,\\
&\mathbb{E}_{a\sim\pi(\cdot\mid s),s'\sim\Pcal(\cdot\mid s,a))}[r(s,a)+\gamma\phi(s')^{\top}\theta\phi(s)-\phi(s)^{\top}\theta\phi(s)-\phi(s)\phi(s)^{\top}\omega]=0.
\end{aligned}\label{eq:TDC}
\end{gather}

A momentum term can be added to the update of the variable $\theta$ to accelerate the convergence of this method. \citet{deb2022gradient} studies this approach and shows that the objective of TDC with momentum can be described as solving a system of three equations. Compared with \eqref{eq:TDC}, the additional time scale is introduced to force the momentum term to decay to zero, which becomes an equation in the lowest level. The second equation in \eqref{eq:TDC} is primarily used to solve for $\omega$ and lies in the middle level, whereas $\theta$ is the variable associated with the highest level equation.
We skip the formulation details here and point interested readers to \citet{deb2022gradient}[Section 4.2]. It is important to note that under a common assumption on the feature vectors (see, for example, \citet{xu2019two}[Assumption 1]), the strong monotonicity condition holds, allowing us to conclude that \texttt{A-MSA} find the (unique) optimal solution with rate $\widetilde{\Ocal}(1/t)$. While TDC with momentum has been shown to converge asymptotically in \citet{deb2022gradient}, its finite-time complexity is unknown from the existing literature. Our paper fills in this gap.

\subsection{Distributionally Reinforcement Q Learning}

Distributionally robust reinforcement learning (DRRL) in general studies finding a policy robust under unknown environmental changes. In this paper, we introduce the subject following the specific formulation in \citet{liang2024single}. 
Consider again the MDP $(\Scal,\Acal,\Pcal,r,\gamma)$ introduced in Section~\ref{sec:applications:TDC}, where $\Pcal$ is the transition kernel that we can sample from during training. The unknown transition kernel in the test environment may deviate from $\Pcal$ but lies within an uncertainty set $\Ucal=\prod_{(s,a)\in\Scal\times\Acal}\Ucal_{s,a}$ with
\[
\Ucal_{s,a}=\{\Pcal'(\cdot\mid s,a)\in\Delta_{\Scal}:D_{\chi^2}(\Pcal'(\cdot\mid s,a),\Pcal(\cdot\mid s,a))\leq\rho\}.
\]
Here $D_{\chi^2}(u,v)$ denotes the $\chi^2$ distance between distributions $u,v\in\Delta_{\Scal}$, and $\rho$ is a radius parameter of the uncertainty set. 
The aim of DRRL is to find the distributionally robust Q function $Q^{\text{rob},\star}\in\mathbb{R}^{|\Scal|\times|\Acal|}$ that performs optimally in the worst case
\begin{align*}
Q^{\text{rob},\star}(s,a)\triangleq \sup_{\pi}\inf_{\Pcal'\in\Ucal}\mathbb{E}_{\pi,\Pcal'}[\sum_{t=0}^{\infty}\gamma^{t}r(s_t,a_t)\mid s_0=s,a_0=a].
\end{align*}

The distributionally robust Q function satisfies the robust Bellman optimality equation
\begin{align}
Q^{\text{rob}, \star}(s, a)=r(s, a)+\gamma \inf_{\Pcal' \in \Ucal} \mathbb{E}_{\Pcal'}[\max_{a \in \Acal} Q^{\text{rob}, \star}(s', a)],\quad\forall s,a.\label{eq:DR_Bellman}
\end{align}
Directly solving $Q^{\text{rob},\star}$ from \eqref{eq:DR_Bellman}, however, is challenging as an exhaustive search in $\Ucal$ is infeasible and we cannot sample from $\Pcal'$. The trick developed in the distributionally robust optimization literature \citep{duchi2021learning} is to use the dual form of the $\chi^2$ distance
\begin{align}
\inf_{D_{\chi^2}\left(\Pcal',\Pcal\right) \leq \rho}\mathbb{E}_{\Pcal}[X]=
\sup_{\eta \in \mathbb{R}}\left\{\eta-\sqrt{1+\rho} \sqrt{\mathbb{E}_{\Pcal}\left[(\eta-X)_{+}^2\right]}\right\}.
\end{align}
We define the shorthand notation $\sigma_{\chi^2}(X, \eta)=\eta-\sqrt{1+\rho} \sqrt{\mathbb{E}_{\Pcal}\left[(\eta-X)_{+}^2\right]}$.
Exploiting this relation, it can be shown that a variant of the Q learning algorithm for the distributionally robust setting aims to find the solution tuple $\{Q^{\text{rob},\star}\in\mathbb{R}^{|\Scal|\times|\Acal|},\vec{\eta}\in\mathbb{R}^{|\Scal|\times|\Acal|},Z_1\in\mathbb{R}^{|\Scal|\times|\Acal|},Z_2\in\mathbb{R}^{|\Scal|\times|\Acal|}\}$ to the system of equations
\begin{align}
&\;Q^{\text{rob},\star}(s,a)=r(s,a)+\gamma (\vec{\eta}(s,a)-\sqrt{1+\rho}\sqrt{Z_2}),\quad\forall s,a,\label{eq:DRRL:obj1}\\
&\;\nabla_{\eta}\sigma_{\chi^2}(\max_{a\in\Acal}Q^{\text{rob},\star}(s,a),\vec{\eta}(s,a))=1-\sqrt{1+\rho}\frac{Z_1(s,a)}{Z_2(s,a)}=0,\quad\forall s,a,\label{eq:DRRL:obj2}\\
&\left\{\begin{aligned}
Z_1(s,a)=\mathbb{E}_{\Pcal}[(\vec{\eta}(s,a)-\max_{a\in\Acal}Q^{\text{rob},\star}(s',a))_+],\\ Z_2(s,a)=\mathbb{E}_{\Pcal}[(\vec{\eta}(s,a)-\max_{a\in\Acal}Q^{\text{rob},\star}(s',a))^2],
\end{aligned}\right.
\quad\forall s,a.\label{eq:DRRL:obj3}
\end{align}
Interested readers can find in \citet{liang2024single} the detailed derivation of the system of equations as well as how a three-time-scale algorithm is developed to solve the system. \citet{liang2024single} further establishes the asymptotic convergence of the algorithm.
The highest level equation is \eqref{eq:DRRL:obj1}, and the associated Bellman operator is strongly monotone under a common contraction assumption in Q learning \citep{chen2022finite}. The lowest level equation is \eqref{eq:DRRL:obj3}, which has an associated operator that is an identity mapping with respect to $Z_1$ and $Z_2$ and is therefore obviously strongly monotone. The issue is that the middle level operator $\nabla_{\eta}\sigma_{\chi^2}$ in \eqref{eq:DRRL:obj2} is in general only (non-strongly) monotone with respect to $\vec{\eta}$. The violation of assumption makes our analysis not immediately applicable. 
However, the techniques developed may help design accelerated three-time-scale algorithms with analysis tailored to the monotone structure of $\nabla_{\eta}\sigma_{\chi^2}$.
Additionally, it is possible that $\nabla_{\eta}\sigma_{\chi^2}$ becomes strongly monotone under additional assumptions on $\Ucal$, in which case we can invoke Theorem~\ref{thm:main} and conclude that the \texttt{A-MSA} algorithm can be applied and find the optimal distributionally robust Q function with rate $\widetilde{\Ocal}(1/t)$. 

\subsection{Actor-Critic Algorithm for Regularized Two-Player Zero-Sum Markov Games}

Consider a two-player zero-sum Markov games described by $(\Scal,\Acal,\Bcal,\Pcal,\gamma,r)$. $\Scal$ is the state space, observed by both players. $\Acal$ and $\Bcal$ are the action spaces of the two players. The transition kernel is $\Pcal:\Scal\times\Acal\times\Bcal\rightarrow\Delta_{\Scal}$, with with $\Pcal(s'\mid s,a,b)$ denoting the probability of the game transitioning from state $s$ to $s'$ when player 1 selects action $a\in\Acal$ and player 2 selects $b\in\Bcal$. The reward function is $r:\Scal\times\Acal\times\Bcal\rightarrow[0,1]$ -- when player 1 and 2 take actions $a$ and $b$ in state $s$, player 1 receives reward $r(s,a,b)$ and player 2 receives $-r(s,a,b)$. 
We denote the two players' policies by $\pi\in\Delta_{\Acal}^{\Scal}$ and $\phi\in\Delta_{\Bcal}^{\Scal}$, with $\pi(a\mid s)$, $\phi(b\mid s)$ representing the probability of selecting action $a$, $b$ in state $s$ according to $\pi$, $\phi$. For a fixed initial state distribution $\rho\in\Delta_{\Scal}$, the expected discounted cumulative reward under entropy regularization (received by player 1) under policy pair $(\pi,\phi)$ is
\begin{align*}
J_{w}(\pi, \phi)=\mathbb{E}_{s_0\sim\rho}[V_{w}^{\pi, \phi}(s_0)],
\end{align*}
where $w\in\mathbb{R}_+$ is a non-negative regularization weight and $V_{w}^{\pi, \phi}(s)$ is the regularized value function
\begin{align*}
V_w^{\pi, \phi}(s)\triangleq\mathbb{E}_{
s_t,a_t,b_t
}\Big[\sum\nolimits_{t=0}^{\infty} \gamma^t \Big(r\left(s_t, a_t, b_t\right)-w\log\pi(a_t\mid s_t)+w\log\phi(b_t\mid s_t)\Big) \mid s_0=s\Big].
\end{align*}
Solving a regularized game means that we want to find a Nash equilibrium policy pair $(\pi_w^{\star},\phi_w^{\star})$ as the solution to 
\begin{align*}
J_{w}(\pi_{w}^{\star},\phi_{w}^{\star})=\max_{\pi\in\Delta_{\Acal}^{\Scal}}\min_{\phi\in\Delta_{\Bcal}^{\Scal}}J_{w}(\pi,\phi)=\min_{\phi\in\Delta_{\Bcal}^{\Scal}}\max_{\pi\in\Delta_{\Acal}^{\Scal}}J_{w}(\pi,\phi).
\end{align*}
It is known from \citet{zeng2022regularized} that such a Nash equilibrium exists and is unique. 

Here we consider the softmax parameterization -- we introduce parameters $\theta\in\mathbb{R}^{|\Scal|\times|\Acal|},\psi\in\mathbb{R}^{|\Scal|\times|\Bcal|}$ that encode the policies $\pi_{\theta},\phi_{\psi}$ according to
\begin{align}
    \pi_{\theta}(a \mid s)=\frac{\exp \left(\theta(s, a)\right)}{\sum_{a' \in \Acal} \exp \left(\theta(s, a')\right)}, \quad \phi_{\psi}(b \mid s)=\frac{\exp \left(\psi(s, b)\right)}{\sum_{b' \in \Acal} \exp \left(\psi(s, b')\right)}.\label{eq:softmax}
\end{align}
Taking a gradient-based approach to the problem, we can express our objective as finding a stationary point where $\nabla_{\theta} J_{w}(\pi_{\theta},\phi_{\psi})=0$ and $\nabla_{\psi} J_{w}(\pi_{\theta},\phi_{\psi})=0$, which can be expanded as below.
\begin{align}
&\mathbb{E}_{s,a,b,s'}\big[(r(s,a)-w\log\pi_{\theta}(a\mid s)+\gamma V_{w}^{\pi_{\theta},\phi_{\psi}}(s')-V_{w}^{\pi_{\theta},\phi_{\psi}}(s))\nabla_{\theta}\log\pi_{\theta}(a\mid s)\big]=0,\label{eq:regularizedMG:eq1}\\
&\mathbb{E}_{s,a,b,s'}\big[(r(s,a)+w\log\phi_{\psi}(b\mid s)+\gamma V_{w}^{\pi_{\theta},\phi_{\psi}}(s')-V_{w}^{\pi_{\theta},\phi_{\psi}}(s))\nabla_{\psi}\log\phi_{\psi}(b\mid s)\big]=0.\label{eq:regularizedMG:eq2}
\end{align}
The value function is not directly known but solvable as the solution to the Bellman equation
\begin{align}
\mathbb{E}_{a,b,s'}\big[r(s,a,b)-w\log\pi_{\theta}(a\mid s)+w\log\phi_{\psi}(b\mid s)+\gamma V_{w}^{\pi_{\theta},\phi_{\psi}}(s')-V_{w}^{\pi_{\theta},\phi_{\psi}}(s)\big]=0,\;\forall s\in\Scal.
\label{eq:regularizedMG:eq3}
\end{align}
We need to solve the system of three equations \eqref{eq:regularizedMG:eq1}-\eqref{eq:regularizedMG:eq3}. The lowest level is \eqref{eq:regularizedMG:eq3} and has corresponds to an operator that can be shown to be strongly monotone. \eqref{eq:regularizedMG:eq1} and \eqref{eq:regularizedMG:eq2} are the highest level and middle level equations, associated with gradient operator $\nabla_{\theta} J_{w}(\pi_{\theta},\phi_{\psi})$ and $\nabla_{\psi} J_{w}(\pi_{\theta},\phi_{\psi})$, which are not strongly monotone. This prevents our analysis from being immediately applicable. However, the operators exhibit an important structure -- they are the gradients of functions satisfying the Polyak-\L ojasiewicz (PL) condition, which makes the operators resemble strong monotone operators in a transformed domain. Exploiting the structure, \citet{zeng2022regularized} shows that a gradient descent ascent algorithm that aims to find the solution to \eqref{eq:regularizedMG:eq1}-\eqref{eq:regularizedMG:eq2} assuming the exact knowledge of $V_w^{\pi_{\theta},\phi_{\psi}}$ converges linearly fast, which is the convergence rate to be expected when the operators $\nabla_{\theta} J_{w}(\pi_{\theta},\phi_{\psi})$ and $\nabla_{\psi} J_{w}(\pi_{\theta},\phi_{\psi})$ are strongly monotone. Combining the techniques in \citet{zeng2022regularized} on leveraging the PL condition and those in this work on acceleration, we believe that the \texttt{A-MSA} algorithm can be shown to find the Nash equilibrium in a regularized two-player zero-sum Markov game with convergence rate $\widetilde{\Ocal}(1/t)$.

\subsection{Actor-Critic Algorithm for Mean Field Games}\label{sec:applications:MFG}

A mean field game (MFG) provides an approximation of an $N$-agent Markov game with homogeneous agents as $N$ approaches infinity. The goal in solving a mean field game is to find a policy optimal in an MDP determined by the induced mean field, where the induced mean field is a distribution over the (infinite) population of agents and a function of the policy itself. In essence, an representative agent in a MFG needs to perform against an infinite population of other agents that adopts its same policy.

We consider MFGs in the stationary infinite-horizon average-reward setting and follow the formulation and notations in \citet{zeng2024single}. We denote the state and action space of the representative agent by $\Scal$ and $\Acal$. The state transition depends not only on the action of the representative agent but also on the aggregate population behavior, which is described by the mean field $u\in\Delta_{\Scal}$. We use $\Pcal:\Scal\times\Acal\times \Delta_{\Scal}\rightarrow\Delta_{\Scal}$ to denote the transition kernel, with $\Pcal(s'\mid s,a,u)$ representing the probability that the state transitions from $s$ to $s'$ when the representative agent takes action $a$ and mean field is $u$. 
Similarly, the reward function $r:\Scal\times\Acal\times\Delta_{\Scal}\rightarrow[0,1]$ depends on the mean field -- the representative agent receives reward $r(s,a,u)$ when it takes action $a$ under mean field $u$ in state $s$. The agent does not observe the mean field, and its policy $\pi$ is a mapping from $\Scal$ to $\Delta_{\Acal}$. 
We denote by $P^{\pi,u}\in\mathbb{R}^{|\Scal|\times|\Scal|}$ the state transition matrix under policy $\pi$ and mean field $u$ such that
\[P^{\pi,u}(s'\mid s)=\sum_{a\in\Acal}\Pcal(s'\mid s,a,u)\pi(a\mid s).\]

When the mean field is $u$ and the agent takes policy $\pi$, the agent can expect to collect the cumulative reward $J(\pi,u)$
\begin{align*}
J(\pi,u)&\triangleq \lim_{T\rightarrow\infty}\frac{1}{T}\textstyle\sum_{t=0}^{T-1}\mathbb{E}_{a_t\sim\pi(\cdot\mid s_t),s_{t+1}\sim\Pcal(\cdot\mid s_t,a_t,u)}[ r(s_t, a_t, u ) \mid s_0].
\end{align*}

The mean field induced by policy $\pi$ (i.e. state visitation distribution over population when every agent takes policy $\pi$) is denoted by $u^{\star}(\pi)$, which satisfies 
\begin{align*}
u^{\star}(\pi)=\lim_{T\rightarrow\infty}\frac{1}{T}\mathbb{E}_{a_t\sim\pi(\cdot\mid s_t),s_{t+1}\sim\Pcal(\cdot\mid s_t,a_t,u^{\star}(\pi))}[e_{s_t}],
\end{align*}
where $e_s\in\mathbb{R}^{|\Scal|}$ is the indicator vector whose entry $s'$ has value 1 if $s'=s$ and 0 otherwise. 

We again consider the softmax function -- a policy parameter $\theta$ encodes the policy $\pi_{\theta}$ according to \eqref{eq:softmax}.
With the derivation presented in \citet{zeng2024single}, it can be shown that finding an equilibrium in an MFG can be formulated as solving a system of three equations \eqref{eq:MFG:eq1}-\eqref{eq:MFG:eq3}
\begin{align}
&\nabla_{\theta} J(\pi_{\theta},u) = \mathbb{E}_{s,a,s'}\big[(r(s,a,u)+V^{\pi_{\theta},u}(s')-V^{\pi_{\theta},u}(s))\nabla_{\theta}\log\pi_{\theta}(a\mid s)\big]=0,\label{eq:MFG:eq1}\\
&u=u^{\star}(\pi_{\theta})=\lim_{T\rightarrow\infty}\frac{1}{T}\mathbb{E}_{a_t\sim\pi_{\theta}(\cdot\mid s_t),s_{t+1}\sim\Pcal(\cdot\mid s_t,a_t,u)}[e_{s_t}],\label{eq:MFG:eq2}\\
&\mathbb{E}_{a,s'}\big[r(s,a,u)-J(\pi_{\theta},u)+V^{\pi_{\theta},u}(s')-V^{\pi_{\theta},u}(s)\big]=0,\;\forall s\in\Scal.\label{eq:MFG:eq3}
\end{align}
in which we introduce the auxiliary variable $V^{\pi_{\theta},u}$ which tracks the (differential) value function under $\pi_{\theta},u$. Among the three equations, \eqref{eq:MFG:eq2} and \eqref{eq:MFG:eq3} are on the middle and lowest levels, and have associated operators satisfying the strong monotonicity. The highest level operator $\nabla_{\theta} J(\pi_{\theta},u)$ is not monotone but a gradient operator. We can extend the analysis in the work leveraging the techniques developed in \citet{zeng2024fast} for gradient operators and show that in this case \texttt{A-MSA} finds a first-order stationary point (but not necessarily a globally or locally optimal solution) of $J$ with rate $\widetilde{\Ocal}(1/\sqrt{t})$ (in the sense of gradient norm squared). This recovers the rate of the state-of-the-art MFG-ASAC algorithm proposed in \citet{zeng2024single} with tailored analysis.

\section{Numerical Simulations}

\begin{figure}[ht]
\centering
\includegraphics[width=\textwidth]{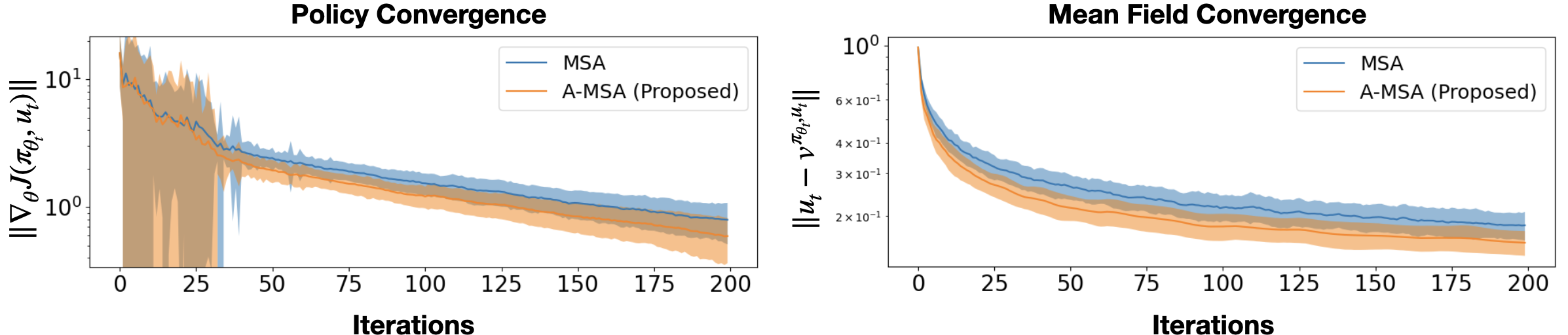}
\vspace{-10pt}
\caption{Algorithm performance in synthetic MFGs averaged over 200 trials. First column measures the sub-optimality in the policy iterate under the latest mean field estimate. Second column shows the convergence of mean field estimate to mean field induced by latest policy.}
\vspace{-5pt}
\label{fig:MFG}
\end{figure}

We apply \texttt{A-MSA} to the MFG policy optimization problem discussed in Section~\ref{sec:applications:MFG}. The environment is a small-scale synthetic MFG with $|\Scal|=30$, $|\Acal|=10$, and a randomly generated transition kernel and reward function. We compare \texttt{A-MSA} against the standard \texttt{MSA} algorithm without averaging and plot the algorithm convergence in Figure~\ref{fig:MFG}, in which
$\theta_t$ and $u_t$ denote the policy and mean field iterates in the $t^{\text{th}}$ iteration. As the norm of the gradient with respect to $\theta_t$, $\|\nabla_{\theta}J(\pi_{\theta_t},u_t)\|$ measures the convergence in the policy given the latest mean field iterate. To quantify the mean field convergence, we consider $\|u_t-v^{\pi_{\theta_t},u_t}\|$, where $v^{\pi,\mu}$ for any $\pi,\mu$ denotes the stationary distribution of states induced by $P^{\pi,\mu}$. The simulation shows that \texttt{A-MSA} has a clear advantage over \texttt{MSA}, though not by an apparent order of magnitude.

\section{Concluding Remarks}
In this work we propose an accelerated multi-time-scale SA algorithm that solves a coupled system of $N$ fixed-point equations with the optimal convergence rate and sample complexity -- under a strong monotonicity assumption, the last iterate of the algorithm converges to the (unique) solution with rate $\widetilde{\Ocal}(1/t)$ for any $N$. This is the first time such convergence rate is achieved under no additional restrictive smoothness assumptions. We show that we can formulate the objectives of a range of problems in RL and multi-agent games as a coupled fixed-point equation system, to which the proposed algorithm can be applied. Our analysis is easily generalizable to cases where the highest level operator is not strongly monotone (but the lower level operators are). An important future direction is to investigate the stability and convergence of multi-time-scale SA when the lower level operators lose strong monotonicity.

\section*{Acknowledgement}
This work was partially supported by the National Science Foundation under awards ECCS-CAREER-2339509 and CCF-2343599.

\section*{Disclaimer}
This paper was prepared for informational purposes in part by
the Artificial Intelligence Research group of JP Morgan Chase \& Co and its affiliates (``JP Morgan''),
and is not a product of the Research Department of JP Morgan.
JP Morgan makes no representation and warranty whatsoever and disclaims all liability,
for the completeness, accuracy or reliability of the information contained herein.
This document is not intended as investment research or investment advice, or a recommendation,
offer or solicitation for the purchase or sale of any security, financial instrument, financial product or service,
or to be used in any way for evaluating the merits of participating in any transaction,
and shall not constitute a solicitation under any jurisdiction or to any person,
if such solicitation under such jurisdiction or to such person would be unlawful.

\medskip

\bibliography{references}

\begin{thebibliography}{44}
\providecommand{\natexlab}[1]{#1}
\providecommand{\url}[1]{\texttt{#1}}
\expandafter\ifx\csname urlstyle\endcsname\relax
  \providecommand{\doi}[1]{doi: #1}\else
  \providecommand{\doi}{doi: \begingroup \urlstyle{rm}\Url}\fi

\bibitem[Benveniste et~al.(2012)Benveniste, M{\'e}tivier, and Priouret]{benveniste2012adaptive}
Albert Benveniste, Michel M{\'e}tivier, and Pierre Priouret.
\newblock \emph{Adaptive algorithms and stochastic approximations}, volume~22.
\newblock Springer Science \& Business Media, 2012.

\bibitem[Bhandari et~al.(2018)Bhandari, Russo, and Singal]{bhandari2018finite}
Jalaj Bhandari, Daniel Russo, and Raghav Singal.
\newblock A finite time analysis of temporal difference learning with linear function approximation.
\newblock In \emph{Conference on learning theory}, pages 1691--1692. PMLR, 2018.

\bibitem[Borkar(1997)]{borkar1997stochastic}
Vivek~S Borkar.
\newblock Stochastic approximation with two time scales.
\newblock \emph{Systems \& Control Letters}, 29\penalty0 (5):\penalty0 291--294, 1997.

\bibitem[Borkar(2008)]{borkar2008stochastic}
Vivek~S Borkar.
\newblock \emph{Stochastic approximation: a dynamical systems viewpoint}.
\newblock Springer, 2008.

\bibitem[Chen et~al.(2022)Chen, Zhang, Doan, Clarke, and Maguluri]{chen2022finite}
Zaiwei Chen, Sheng Zhang, Thinh~T Doan, John-Paul Clarke, and Siva~Theja Maguluri.
\newblock Finite-sample analysis of nonlinear stochastic approximation with applications in reinforcement learning.
\newblock \emph{Automatica}, 146:\penalty0 110623, 2022.

\bibitem[Dalal et~al.(2018)Dalal, Thoppe, Sz{\"o}r{\'e}nyi, and Mannor]{dalal2018finite}
Gal Dalal, Gugan Thoppe, Bal{\'a}zs Sz{\"o}r{\'e}nyi, and Shie Mannor.
\newblock Finite sample analysis of two-timescale stochastic approximation with applications to reinforcement learning.
\newblock In \emph{Conference On Learning Theory}, pages 1199--1233. PMLR, 2018.

\bibitem[Daskalakis et~al.(2020)Daskalakis, Foster, and Golowich]{daskalakis2020independent}
Constantinos Daskalakis, Dylan~J Foster, and Noah Golowich.
\newblock Independent policy gradient methods for competitive reinforcement learning.
\newblock \emph{Advances in neural information processing systems}, 33:\penalty0 5527--5540, 2020.

\bibitem[Deb and Bhatnagar(2021)]{deb2021n}
Rohan Deb and Shalabh Bhatnagar.
\newblock {$N$}-timescale stochastic approximation: Stability and convergence.
\newblock \emph{arXiv preprint arXiv:2112.03515}, 2021.

\bibitem[Deb and Bhatnagar(2022)]{deb2022gradient}
Rohan Deb and Shalabh Bhatnagar.
\newblock Gradient temporal difference with momentum: Stability and convergence.
\newblock In \emph{Proceedings of the AAAI Conference on Artificial Intelligence}, volume~36, pages 6488--6496, 2022.

\bibitem[Doan(2022)]{doan2022nonlinear}
Thinh~T Doan.
\newblock Nonlinear two-time-scale stochastic approximation: Convergence and finite-time performance.
\newblock \emph{IEEE Transactions on Automatic Control}, 68\penalty0 (8):\penalty0 4695--4705, 2022.

\bibitem[Doan(2024)]{doan2024fast}
Thinh~T Doan.
\newblock Fast nonlinear two-time-scale stochastic approximation: Achieving $\mathcal{O}(1/k)$ finite-sample complexity.
\newblock \emph{arXiv preprint arXiv:2401.12764}, 2024.

\bibitem[Doan and Romberg(2019)]{doan2019linear}
Thinh~T Doan and Justin Romberg.
\newblock Linear two-time-scale stochastic approximation a finite-time analysis.
\newblock In \emph{2019 57th Annual Allerton Conference on Communication, Control, and Computing (Allerton)}, pages 399--406. IEEE, 2019.

\bibitem[Duchi and Namkoong(2021)]{duchi2021learning}
John~C Duchi and Hongseok Namkoong.
\newblock Learning models with uniform performance via distributionally robust optimization.
\newblock \emph{The Annals of Statistics}, 49\penalty0 (3):\penalty0 1378--1406, 2021.

\bibitem[Fort(2015)]{fort2015central}
Gersende Fort.
\newblock Central limit theorems for stochastic approximation with controlled markov chain dynamics.
\newblock \emph{ESAIM: Probability and Statistics}, 19:\penalty0 60--80, 2015.

\bibitem[Gupta et~al.(2019)Gupta, Srikant, and Ying]{gupta2019finite}
Harsh Gupta, Rayadurgam Srikant, and Lei Ying.
\newblock Finite-time performance bounds and adaptive learning rate selection for two time-scale reinforcement learning.
\newblock \emph{Advances in Neural Information Processing Systems}, 32, 2019.

\bibitem[Han et~al.(2024)Han, Li, and Zhang]{han2024finite}
Yuze Han, Xiang Li, and Zhihua Zhang.
\newblock Finite-time decoupled convergence in nonlinear two-time-scale stochastic approximation.
\newblock \emph{arXiv preprint arXiv:2401.03893}, 2024.

\bibitem[Haque et~al.(2023)Haque, Khodadadian, and Maguluri]{haque2023tight}
Shaan~Ul Haque, Sajad Khodadadian, and Siva~Theja Maguluri.
\newblock Tight finite time bounds of two-time-scale linear stochastic approximation with markovian noise.
\newblock \emph{arXiv preprint arXiv:2401.00364}, 2023.

\bibitem[Hong et~al.(2023)Hong, Wai, Wang, and Yang]{hong2023two}
Mingyi Hong, Hoi-To Wai, Zhaoran Wang, and Zhuoran Yang.
\newblock A two-timescale stochastic algorithm framework for bilevel optimization: Complexity analysis and application to actor-critic.
\newblock \emph{SIAM Journal on Optimization}, 33\penalty0 (1):\penalty0 147--180, 2023.

\bibitem[Kaledin et~al.(2020)Kaledin, Moulines, Naumov, Tadic, and Wai]{kaledin2020finite}
Maxim Kaledin, Eric Moulines, Alexey Naumov, Vladislav Tadic, and Hoi-To Wai.
\newblock Finite time analysis of linear two-timescale stochastic approximation with markovian noise.
\newblock In \emph{Conference on Learning Theory}, pages 2144--2203. PMLR, 2020.

\bibitem[Konda and Tsitsiklis(1999)]{konda1999actorb}
Vijay Konda and John Tsitsiklis.
\newblock Actor-critic algorithms.
\newblock \emph{Advances in neural information processing systems}, 12, 1999.

\bibitem[Konda and Tsitsiklis(2004)]{konda2004convergence}
Vijay~R Konda and John~N Tsitsiklis.
\newblock Cconvergence rate of linear two-time-scale stochastic approximation.
\newblock \emph{The Annals of Applied Probability}, 14\penalty0 (2):\penalty0 796--819, 2004.

\bibitem[Konda and Borkar(1999)]{konda1999actora}
Vijaymohan~R Konda and Vivek~S Borkar.
\newblock Actor-critic--type learning algorithms for markov decision processes.
\newblock \emph{SIAM Journal on control and Optimization}, 38\penalty0 (1):\penalty0 94--123, 1999.

\bibitem[Lakshminarayanan and Szepesv{\'a}ri(2017)]{lakshminarayanan2017linear}
Chandrashekar Lakshminarayanan and Csaba Szepesv{\'a}ri.
\newblock Linear stochastic approximation: Constant step-size and iterate averaging.
\newblock \emph{arXiv preprint arXiv:1709.04073}, 2017.

\bibitem[Li et~al.(2023{\natexlab{a}})Li, Liang, and Zhang]{li2023online}
Xiang Li, Jiadong Liang, and Zhihua Zhang.
\newblock Online statistical inference for nonlinear stochastic approximation with markovian data.
\newblock \emph{arXiv preprint arXiv:2302.07690}, 2023{\natexlab{a}}.

\bibitem[Li et~al.(2023{\natexlab{b}})Li, Yang, Liang, Zhang, and Jordan]{li2023statistical}
Xiang Li, Wenhao Yang, Jiadong Liang, Zhihua Zhang, and Michael~I Jordan.
\newblock A statistical analysis of polyak-ruppert averaged q-learning.
\newblock In \emph{International Conference on Artificial Intelligence and Statistics}, pages 2207--2261. PMLR, 2023{\natexlab{b}}.

\bibitem[Liang et~al.(2024)Liang, Ma, Blanchet, Yang, Zhang, and Zhou]{liang2024single}
Zhipeng Liang, Xiaoteng Ma, Jose Blanchet, Jun Yang, Jiheng Zhang, and Zhengyuan Zhou.
\newblock Single-trajectory distributionally robust reinforcement learning.
\newblock In \emph{Proceedings of the 41st International Conference on Machine Learning}, volume 235, pages 29644--29666. PMLR, 2024.

\bibitem[Meerkov(1972)]{Meerkov1972a}
S.~M. Meerkov.
\newblock Simplified description of slow markov walks: part i.
\newblock \emph{Automation and Remote Control}, 33:\penalty0 404--414, 1972.

\bibitem[Mitra(2024)]{mitra2024simple}
Aritra Mitra.
\newblock A simple finite-time analysis of td learning with linear function approximation.
\newblock \emph{arXiv preprint arXiv:2403.02476}, 2024.

\bibitem[Mokkadem and Pelletier(2006)]{mokkadem2006convergence}
Abdelkader Mokkadem and Mariane Pelletier.
\newblock Convergence rate and averaging of nonlinear two-time-scale stochastic approximation algorithms.
\newblock \emph{Annals of Applied Probability}, 16\penalty0 (3):\penalty0 1671--1702, 2006.

\bibitem[Moulines and Bach(2011)]{moulines2011non}
Eric Moulines and Francis Bach.
\newblock Non-asymptotic analysis of stochastic approximation algorithms for machine learning.
\newblock \emph{Advances in neural information processing systems}, 24, 2011.

\bibitem[Robbins and Monro(1951)]{robbins1951stochastic}
Herbert Robbins and Sutton Monro.
\newblock A stochastic approximation method.
\newblock \emph{The annals of mathematical statistics}, pages 400--407, 1951.

\bibitem[Sayin et~al.(2021)Sayin, Zhang, Leslie, Basar, and Ozdaglar]{sayin2021decentralized}
Muhammed Sayin, Kaiqing Zhang, David Leslie, Tamer Basar, and Asuman Ozdaglar.
\newblock Decentralized q-learning in zero-sum markov games.
\newblock \emph{Advances in Neural Information Processing Systems}, 34:\penalty0 18320--18334, 2021.

\bibitem[Shen and Chen(2022)]{shen2022single}
Han Shen and Tianyi Chen.
\newblock A single-timescale analysis for stochastic approximation with multiple coupled sequences.
\newblock \emph{Advances in Neural Information Processing Systems}, 35:\penalty0 17415--17429, 2022.

\bibitem[Srikant(2024)]{srikant2024rates}
R~Srikant.
\newblock Rates of convergence in the central limit theorem for markov chains, with an application to td learning.
\newblock \emph{arXiv preprint arXiv:2401.15719}, 2024.

\bibitem[Srikant and Ying(2019)]{srikant2019finite}
Rayadurgam Srikant and Lei Ying.
\newblock Finite-time error bounds for linear stochastic approximation andtd learning.
\newblock In \emph{Conference on Learning Theory}, pages 2803--2830. PMLR, 2019.

\bibitem[Sutton et~al.(2008)Sutton, Maei, and Szepesv{\'a}ri]{sutton2008convergent}
Richard~S Sutton, Hamid Maei, and Csaba Szepesv{\'a}ri.
\newblock A convergent {$O(n)$} temporal-difference algorithm for off-policy learning with linear function approximation.
\newblock \emph{Advances in neural information processing systems}, 21, 2008.

\bibitem[Sutton et~al.(2009)Sutton, Maei, Precup, Bhatnagar, Silver, Szepesv{\'a}ri, and Wiewiora]{sutton2009fast}
Richard~S Sutton, Hamid~Reza Maei, Doina Precup, Shalabh Bhatnagar, David Silver, Csaba Szepesv{\'a}ri, and Eric Wiewiora.
\newblock Fast gradient-descent methods for temporal-difference learning with linear function approximation.
\newblock In \emph{Proceedings of the 26th annual international conference on machine learning}, pages 993--1000, 2009.

\bibitem[Wang et~al.(2021)Wang, Zou, and Zhou]{wang2021non}
Yue Wang, Shaofeng Zou, and Yi~Zhou.
\newblock Non-asymptotic analysis for two time-scale tdc with general smooth function approximation.
\newblock \emph{Advances in neural information processing systems}, 34:\penalty0 9747--9758, 2021.

\bibitem[Wu et~al.(2020)Wu, Zhang, Xu, and Gu]{wu2020finite}
Yue~Frank Wu, Weitong Zhang, Pan Xu, and Quanquan Gu.
\newblock A finite-time analysis of two time-scale actor-critic methods.
\newblock \emph{Advances in Neural Information Processing Systems}, 33:\penalty0 17617--17628, 2020.

\bibitem[Xu et~al.(2019)Xu, Zou, and Liang]{xu2019two}
Tengyu Xu, Shaofeng Zou, and Yingbin Liang.
\newblock Two time-scale off-policy td learning: Non-asymptotic analysis over markovian samples.
\newblock \emph{Advances in neural information processing systems}, 32, 2019.

\bibitem[Zeng and Doan(2024)]{zeng2024fast}
Sihan Zeng and Thinh Doan.
\newblock Fast two-time-scale stochastic gradient method with applications in reinforcement learning.
\newblock In \emph{The Thirty Seventh Annual Conference on Learning Theory}, pages 5166--5212. PMLR, 2024.

\bibitem[Zeng et~al.(2022)Zeng, Doan, and Romberg]{zeng2022regularized}
Sihan Zeng, Thinh Doan, and Justin Romberg.
\newblock Regularized gradient descent ascent for two-player zero-sum markov games.
\newblock \emph{Advances in Neural Information Processing Systems}, 35:\penalty0 34546--34558, 2022.

\bibitem[Zeng et~al.(2024{\natexlab{a}})Zeng, Bhatt, Koppel, and Ganesh]{zeng2024single}
Sihan Zeng, Sujay Bhatt, Alec Koppel, and Sumitra Ganesh.
\newblock A single-loop finite-time convergent policy optimization algorithm for mean field games (and average-reward markov decision processes).
\newblock \emph{arXiv preprint arXiv:2408.04780}, 2024{\natexlab{a}}.

\bibitem[Zeng et~al.(2024{\natexlab{b}})Zeng, Doan, and Romberg]{zeng2024two}
Sihan Zeng, Thinh~T Doan, and Justin Romberg.
\newblock A two-time-scale stochastic optimization framework with applications in control and reinforcement learning.
\newblock \emph{SIAM Journal on Optimization}, 34\penalty0 (1):\penalty0 946--976, 2024{\natexlab{b}}.

\end{thebibliography}

\clearpage
\vbox{%
\hsize\textwidth
\linewidth\hsize
\hrule height 1pt
\vskip 0.25in
\centering
{\Large\bf{Accelerated Multi-Time-Scale Stochastic Approximation: Optimal Complexity and Applications in Reinforcement Learning and Multi-Agent Games \\ 
Appendix} \par}
  \vskip 0.25in
\hrule height 1pt
\vskip 0.09in
}

\tableofcontents

\appendix

\section{Proof of Technical Lemmas}\label{sec:proof_lemmas}

\subsection{Proof of Lemma~\ref{lem:bound_Deltaf}}

\begin{lem}\label{lem:Markovian_crossterm}
For any $t\geq\tau_t$ we have
\begin{align*}
\mathbb{E}[\langle\Delta f_i^{[t]},F_i(\vtheta^{[t]},X^{[t]})-\widebar{F}_i(\vtheta^{[t]})\rangle]\leq N^2 C_{\tau}\tau_t^2\lambda^{[t-\tau_t]}\Big(\sum_{j=1}^{N}\|x_j^{[t]}\|^2+\sum_{j=1}^{N}\|\Delta f_j^{[t]}\|^2+1\Big),
\end{align*}
where $C_{\tau}$ is a positive, bounded scalar depending only on the structural constants including $L$, $B$, $\rho$, and $m$.
\end{lem}

We derive a frequently used inequality in the proofs of Lemma~\ref{lem:bound_Deltaf} and \ref{lem:Markovian_crossterm} as a result of \eqref{eq:bounded_F}. For any $X$, we have
\begin{align}
\|F_i(\vtheta^{[t]},X)\|&\leq L(\sum_{j=1}^N\|\theta_j^{[t]}\|+1)\notag\\
&\leq L+L\sum_{j=1}^N\|\theta_j^{[t]}-y_j(\vtheta_{1:j-1}^{[t]})\|+L\sum_{j=1}^N\|y_j(\vtheta_{1:j-1}^{[t]})\|\notag\\
&\leq (NB+1)L+L\sum_{j=1}^{N}\|x_j^{[t]}\|,
\label{lem:bound_Deltaf:boundedF}
\end{align}
where the last inequality follows from the boundedness condition in Assumption~\ref{assump:Lipschitz} and the definition of $x_j^{[t]}$ in \eqref{eq:residuals}.

By the update rule in \eqref{eq:update_auxiliary}, we have
\begin{align*}
    \Delta f_i^{[t+1]}&=f_i^{[t+1]}-\widebar{F}_i(\vtheta^{[t+1]})\notag\\
    &=(1-\lambda^{[t]})f_i^{[t]}+\lambda^{[t]}F_i(\vtheta^{[t]},X^{[t]})-\widebar{F}_i(\vtheta^{[t+1]})\notag\\
    &=(1-\lambda^{[t]})f_i^{[t]}+\lambda^{[t]}\widebar{F}_i(\vtheta^{[t]})-\widebar{F}_i(\vtheta^{[t+1]})+\lambda^{[t]}\left(F_i(\vtheta^{[t]},X^{[t]})-\widebar{F}_i(\vtheta^{[t]})\right)\notag\\
    &=(1-\lambda^{[t]})\Delta f_i^{[t]}+\left(\widebar{F}_i(\vtheta^{[t]})-\widebar{F}_i(\vtheta^{[t+1]})\right)+\lambda^{[t]}\left(F_i(\vtheta^{[t]},X^{[t]})-\widebar{F}_i(\vtheta^{[t]})\right).
\end{align*}

This leads to
\begin{align}
&\|\Delta f_i^{[t+1]}\|^2\notag\\
&=(1-\lambda^{[t]})^2\|\Delta f_i^{[t]}\|^2+\|\left(\widebar{F}_i(\vtheta^{[t]})-\widebar{F}_i(\vtheta^{[t+1]})\right)+\lambda^{[t]}\left(F_i(\vtheta^{[t]},X^{[t]})-\widebar{F}_i(\vtheta^{[t]})\right)\|^2\notag\\
&\hspace{20pt}+2(1-\lambda^{[t]})\langle\Delta f_i^{[t]},\widebar{F}_i(\vtheta^{[t]})-\widebar{F}_i(\vtheta^{[t+1]})\rangle\notag\\
&\hspace{20pt}+2\lambda^{[t]}(1-\lambda^{[t]})\langle\Delta f_i^{[t]},F_i(\vtheta^{[t]},X^{[t]})-\widebar{F}_i(\vtheta^{[t]})\rangle\notag\\
&\leq(1-\lambda^{[t]})^2\|\Delta f_i^{[t]}\|^2+2\|\widebar{F}_i(\vtheta^{[t]})-\widebar{F}_i(\vtheta^{[t+1]})\|^2+2(\lambda^{[t]})^2\|\widebar{F}_i(\vtheta^{[t]},X^{[t]})-\widebar{F}_i(\vtheta^{[t]})\|^2\notag\\
&\hspace{20pt}+\frac{\lambda^{[t]}}{2}\|\Delta f_i^{[t]}\|^2+\frac{2}{\lambda^{[t]}}\|\widebar{F}_i(\vtheta^{[t]})-\widebar{F}_i(\vtheta^{[t+1]})\|^2\notag\\
&\hspace{20pt}+2\lambda^{[t]}(1-\lambda^{[t]})\langle\Delta f_i^{[t]},F_i(\vtheta^{[t]},X^{[t]})-\widebar{F}_i(\vtheta^{[t]})\rangle\notag\\
&\leq(1-\lambda^{[t]})\|\Delta f_i^{[t]}\|^2-\frac{\lambda^{[t]}}{4}\|\Delta f_i^{[t]}\|^2+2(\lambda^{[t]})^2\|F_i(\vtheta^{[t]},X^{[t]})-\widebar{F}_i(\vtheta^{[t]})\|^2\notag\\
&\hspace{20pt}+\frac{4}{\lambda^{[t]}}\|\widebar{F}_i(\vtheta^{[t]})-\widebar{F}_i(\vtheta^{[t+1]})\|^2+2\lambda^{[t]}(1-\lambda^{[t]})\langle\Delta f_i^{[t]},F_i(\vtheta^{[t]},X^{[t]})-\widebar{F}_i(\vtheta^{[t]})\rangle,
\label{lem:bound_Deltaf:proof_eq1}
\end{align}
where the last inequality follows from $\lambda^{[t]}\leq1/4\leq 1$. 

Note that the third term on the right hand side of \eqref{lem:bound_Deltaf:proof_eq1} is bounded due to \eqref{lem:bound_Deltaf:boundedF}
\begin{align}
\|F_i(\vtheta^{[t]},X^{[t]})-\widebar{F}_i(\vtheta^{[t]})\|^2 &\leq  2\|F_i(\vtheta^{[t]},X^{[t]})\|+2\|\widebar{F}_i(\vtheta^{[t]})\|^2\notag\\
&\leq 4\Big((NB+1)L+L\sum_{j=1}^{N}\|x_j^{[t]}\|\Big)^2\notag\\
&\leq 8L^2\big(N\sum_{j=1}^N\|x_j^{[t]}\|^2+(NB+1)^2\big).
\label{lem:bound_Deltaf:proof_eq1.1}
\end{align}

For the fourth term of \eqref{lem:bound_Deltaf:proof_eq1}, we have from Lemma~\ref{lem:Lipschitz} and the Lipschitz continuity of $\widebar{F}_i$ 
\begin{align}
\mathbb{E}[\|\widebar{F}_i(\vtheta^{[t]})-\widebar{F}_i(\vtheta^{[t+1]})\|^2]&\leq L^2\sum_{j=1}^{N}\mathbb{E}[\|\theta_j^{[t]}-\theta_j^{[t+1]}\|^2]\notag\\
&\leq L^2(\alpha_i^{[t]})^2\sum_{j=1}^{N}\mathbb{E}[\left(\|\Delta f_j^{[t]}\|+NL^2\sum_{k=j}^{N}\|x_k^{[t]}\|\right)^2]\notag\\
&\leq 2L^2(\alpha_i^{[t]})^2\sum_{j=1}^{N}\mathbb{E}[\|\Delta f_j^{[t]}\|^2]+2NL^4(\alpha_i^{[t]})^2\sum_{j=1}^{N}\sum_{k=j}^{N}\mathbb{E}[\|x_k^{[t]}\|^2]\notag\\
&\leq 2L^2(\alpha_i^{[t]})^2\sum_{j=1}^{N}\mathbb{E}[\|\Delta f_j^{[t]}\|^2]+2N^2L^4(\alpha_i^{[t]})^2\sum_{j=1}^{N}\mathbb{E}[\|x_j^{[t]}\|^2].
\end{align}

The last term of \eqref{lem:bound_Deltaf:proof_eq1} is bound in expectation due to Lemma~\ref{lem:Markovian_crossterm}
\begin{align}
\mathbb{E}[\langle\Delta f_i^{[t]},F_i(\vtheta^{[t]},X^{[t]})-\widebar{F}_i(\vtheta^{[t]})\rangle]\leq N^2 C_{\tau}\tau_t^2\lambda^{[t-\tau_t]}\Big(\sum_{j=1}^{N}\|x_j^{[t]}\|^2+\sum_{j=1}^{N}\|\Delta f_j^{[t]}\|^2+1\Big).
\label{lem:bound_Deltaf:proof_eq1.2}
\end{align}

Combining \eqref{lem:bound_Deltaf:proof_eq1.1}-\eqref{lem:bound_Deltaf:proof_eq1.2} with
\eqref{lem:bound_Deltaf:proof_eq1},
\begin{align}
    &\mathbb{E}[\|\Delta f_i^{[t+1]}\|^2]\notag\\
    &\leq (1-\lambda^{[t]})\mathbb{E}[\|\Delta f_i^{[t]}\|^2]-\frac{\lambda^{[t]}}{4}\mathbb{E}[\|\Delta f_i^{[t]}\|^2]+2(\lambda^{[t]})^2\mathbb{E}[\|F_i(\vtheta^{[t]},X^{[t]})-\widebar{F}_i(\vtheta^{[t]})\|^2]\notag\\
    &\hspace{20pt}+\frac{4}{\lambda^{[t]}}\mathbb{E}[\|\widebar{F}_i(\vtheta^{[t]})-\widebar{F}_i(\vtheta^{[t+1]})\|^2]+2\lambda^{[t]}(1-\lambda^{[t]})\mathbb{E}[\langle\Delta f_i^{[t]},F_i(\vtheta^{[t]},X^{[t]})-\widebar{F}_i(\vtheta^{[t]})\rangle]\notag\\
    &\leq (1-\lambda^{[t]})\mathbb{E}[\|\Delta f_i^{[t]}\|^2]-\frac{\lambda^{[t]}}{4}\mathbb{E}[\|\Delta f_i^{[t]}\|^2]+8L^2(\lambda^{[t]})^2\big(N\sum_{j=1}^N\|x_j^{[t]}\|^2+(NB+1)^2\big)\notag\\
    &\hspace{20pt}+\frac{8L^2(\alpha_i^{[t]})^2}{\lambda^{[t]}}\sum_{j=1}^{N}\mathbb{E}[\|\Delta f_j^{[t]}\|^2] + \frac{8N^2 L^4(\alpha_i^{[t]})^2}{\lambda^{[t]}}\sum_{j=1}^{N}\mathbb{E}[\|x_j^{[t]}\|^2]\notag\\
    &\hspace{20pt}+N^2 C_{\tau}\tau_t^2\lambda^{[t]}\lambda^{[t-\tau_t]}(\sum_{j=1}^{N}\|x_j^{[t]}\|^2+\sum_{j=1}^{N}\|\Delta f_j^{[t]}\|^2+1)\notag\\
    &\leq (1-\lambda^{[t]})\mathbb{E}[\|\Delta f_i^{[t]}\|^2]-\frac{\lambda^{[t]}}{4}\mathbb{E}[\|\Delta f_i^{[t]}\|^2]+DN^2\tau_t^2\lambda^{[t]}\lambda^{[t-\tau_t]}\notag\\
    &\hspace{20pt}+D N^2 \tau_t^2  (\lambda^{[t]}\lambda^{[t-\tau_t]}+\frac{(\alpha_i^{[t]})^2}{\lambda^{[t]}})\big(\sum_{j=1}^{N}\|x_j^{[t]}\|^2+\sum_{j=1}^{N}\|\Delta f_j^{[t]}\|^2\big),
\end{align}
where $D$ is a finite constant depending only on $L$, $B$, and $C_{\tau}$.

\qed

\subsection{Proof of Lemma~\ref{lem:bound_x}}\label{sec:proof_lem:bound_x}

By the update rule in \eqref{alg:update_decision}, we have
\begin{align*}
x_i^{[t+1]}&=\theta_{i}^{[t+1]}-y_i(\vtheta_{1:i-1}^{[t+1]})\notag\\
&=\theta_{i}^{[t]}-\alpha_i^{[t]}f_i^{[t]}-y_i(\vtheta_{1:i-1}^{[t+1]})\notag\\
&=x_i^{[t]}-\alpha_i^{[t]}f_i^{[t]}+\left(y_i(\vtheta_{1:i-1}^{[t]})-y_i(\vtheta_{1:i-1}^{[t+1]})\right)\notag\\
&=x_i^{[t]}-\alpha_i^{[t]} \widebar{F}_i(\vtheta^{[t]})-\alpha_i^{[t]}\Delta f_i^{[t]}+\left(y_i(\vtheta_{1:i-1}^{[t]})-y_i(\vtheta_{1:i-1}^{[t+1]})\right).
\end{align*}

Taking the norm leads to
\begin{align}
&\|x_i^{[t+1]}\|^2 \notag\\
& \leq \underbrace{\|x_i^{[t]}-\alpha_i^{[t]}\widebar{F}_i(\vtheta^{[t]})\|^2}_{T_1} + (\alpha_i^{[t]})^2\|\Delta f_i^{[t]}\|^2 + \underbrace{\|y_i(\vtheta_{1:i-1}^{[t]})-y_i(\vtheta_{1:i-1}^{[t+1]})\|^2}_{T_2}\notag\\
&\hspace{20pt}+\underbrace{\alpha_i^{[t]}\langle x_i^{[t]}-\alpha_i^{[t]}\widebar{F}_i(\vtheta^{[t]}),\Delta f_i^{[t]}\rangle}_{T_3}+\underbrace{\langle x_i^{[t]}-\alpha_i^{[t]}\widebar{F}_i(\vtheta^{[t]}),y_i(\vtheta_{1:i-1}^{[t]})-y_i(\vtheta_{1:i-1}^{[t+1]})\rangle}_{T_4}\notag\\
&\hspace{20pt}+\underbrace{\alpha_i^{[t]}\langle\Delta f_i^{[t]},y_i(\vtheta_{1:i-1}^{[t]})-y_i(\vtheta_{1:i-1}^{[t+1]})\rangle}_{T_5}.\label{lem:bound_x:proof_eq1}
\end{align}

We bound each term of \eqref{lem:bound_x:proof_eq1} separately. By the definition of $y_k$, we have
\begin{align*}
\widebar{F}_i(\vtheta_{1:i-1}^{[t]},\ybf_{i:N}(\vtheta_{1:i-1}^{[t]}))=0.
\end{align*}
This implies
\begin{align}
&T_1\notag\\
&=\|x_i^{[t]}\|^2-2\alpha_i^{[t]}\langle x_i^{[t]},F_i(\vtheta^{[t]})\rangle+(\alpha_i^{[t]})^2\|F_i(\vtheta^{[t]})\|^2\notag\\
&=\|x_i^{[t]}\|^2-2\alpha_i^{[t]}\langle x_i^{[t]},\widebar{F}_i(\vtheta_{1:i}^{[t]},\ybf_{i+1:N}(\vtheta_{1:i}^{[t]}))-\widebar{F}_i(\vtheta_{1:i-1}^{[t]},\ybf_{i:N}(\vtheta_{1:i-1}^{[t]}))\rangle\notag\\
&\hspace{20pt}-2\alpha_i^{[t]}\langle x_i^{[t]}, \widebar{F}_i(\vtheta^{[t]})-\widebar{F}_i(\vtheta_{1:i}^{[t]},\ybf_{i+1:N}(\vtheta_{1:i}^{[t]}))\rangle+(\alpha_i^{[t]})^2\|\widebar{F}_i(\vtheta^{[t]})\|^2\notag\\
&=\|x_i^{[t]}\|^2-2\alpha_i^{[t]}\Big\langle x_i^{[t]},\widebar{F}_i(\vtheta_{1:i-1}^{[t]},\theta_i^{[t]},\ybf_{i+1:N}(\vtheta_{1:i-1}^{[t]},\theta_i^{[t]}))\notag\\
&\hspace{100pt}-\widebar{F}_i\Big(\vtheta_{1:i-1}^{[t]},y_i(\vtheta_{1:i-1}^{[t]}),\ybf_{i+1:N}(\vtheta_{1:i-1}^{[t]}, y_i(\vtheta_{1:i-1}^{[t]}))\Big)\Big\rangle\notag\\
&\hspace{20pt}-2\alpha_i^{[t]}\langle x_i^{[t]}, \widebar{F}_i(\vtheta^{[t]})-\widebar{F}_i(\vtheta_{1:i}^{[t]},\ybf_{i+1:N}(\vtheta_{1:i}^{[t]}))\rangle\notag\\
&\hspace{20pt}+(\alpha_i^{[t]})^2\|\widebar{F}_i(\vtheta^{[t]})-\widebar{F}_i(\vtheta_{1:i-1}^{[t]},\ybf_{i:N}(\vtheta_{1:i-1}^{[t]}))\|^2\notag\\
&\leq\|x_i^{[t]}\|^2-2\delta\alpha_i^{[t]}\|x_i^{[t]}\|^2+2\alpha_i^{[t]}\|x_i^{[t]}\|\|\widebar{F}_i(\vtheta^{[t]})-\widebar{F}_i(\vtheta_{1:i}^{[t]},\ybf_{i+1:N}(\vtheta_{1:i}^{[t]}))\|\notag\\
&\hspace{20pt}+(\alpha_i^{[t]})^2\|\widebar{F}_i(\vtheta^{[t]})-\widebar{F}_i(\vtheta_{1:i-1}^{[t]},\ybf_{i:N}(\vtheta_{1:i-1}^{[t]}))\|^2\notag\\
&\leq\|x_i^{[t]}\|^2-2\delta\alpha_i^{[t]}\|x_i^{[t]}\|^2+\delta\alpha_i^{[t]}\|x_i^{[t]}\|^2+\frac{4N^3L^6\alpha_i^{[t]}}{\delta}\sum_{j=i+1}^{N}\|x_{j}^{[t]}\|^2\notag\\
&\hspace{20pt}+4N^3 L^6(\alpha_i^{[t]})^2\sum_{j=i}^{N}\|x_j^{[t]}\|^2\notag\\
&= \|x_i^{[t]}\|^2-\left(\delta\alpha_i^{[t]}-4N^3L^6(\alpha_i^{[t]})^2\right)\|x_i^{[t]}\|^2+\left(\frac{4N^3L^6\alpha_i^{[t]}}{\delta}+4N^3 L^6(\alpha_i^{[t]})^2\right)\sum_{j=i+1}^{N}\|x_j^{[t]}\|^2\notag\\
&\leq \|x_i^{[t]}\|^2-\frac{\delta\alpha_i^{[t]}}{2}\|x_i^{[t]}\|^2+\frac{8N^3 L^6\alpha_i^{[t]}}{\delta}\sum_{j=i+1}^{N}\|x_j^{[t]}\|^2.
\label{lem:bound_x:proof_eq2}
\end{align}
The third equation of \eqref{lem:bound_x:proof_eq2} is a result of the identity that for all $i,j$
\begin{align}
y_j(\vtheta_{1:i-1})=y_j(\vtheta_{1:i-1}, \ybf_{i:j-1}(\vtheta_{1:i-1})).\label{eq:y_transformation}
\end{align}
The first inequality of \eqref{lem:bound_x:proof_eq2} follows from Assumption~\ref{assump:stronglymonotone}, and the final inequality follows from the step size condition $\alpha_i^{[t]}\leq \frac{\delta}{8N^3 L^6}$ and $\alpha_i^{[t]}\leq \frac{1}{\delta}$.
The second inequality of \eqref{lem:bound_x:proof_eq2} uses the bound 
\[\|\widebar{F}_i(\vtheta^{[t]})-\widebar{F}_i(\vtheta_{1:i}^{[t]},\ybf_{i+1:N}(\vtheta_{1:i}^{[t]}))\|\leq2NL^3\sum_{j=i+1}^{N}\|x_j^{[t]}\|,\] 
which we now justify
\begin{align}
&\|\widebar{F}_i(\vtheta^{[t]})-\widebar{F}_i(\vtheta_{1:i}^{[t]},\ybf_{i+1:N}(\vtheta_{1:i}^{[t]}))\|\notag\\
&\leq L\sum_{j=i+1}^{N}\|\theta_j^{[t]}-y_k(\vtheta_{1:i}^{[t]})\|\notag\\
&\leq L\sum_{j=i+1}^{N}\|\theta_j^{[t]}-y_j(\vtheta_{1:j-1}^{[t]})\|+L\sum_{j=i+1}^{N}\|y_j(\vtheta_{1:j-1}^{[t]})-y_j(\vtheta_{1:i}^{[t]})\|\notag\\
&= L\sum_{j=i+1}^{N}\|x_j^{[t]}\|+L\sum_{j=i+1}^{N}\|y_j(\vtheta_{1:j-1}^{[t]})-y_j(\vtheta_{1:i}^{[t]},\ybf_{i:j-1}(\vtheta_{1:i}^{[t]}))\|\notag\\
&\leq L\sum_{j=i+1}^{N}\|x_j^{[t]}\|+L^3\sum_{j=i+1}^{N}\sum_{k=i+1}^{j-1}\|x_k^{[t]}\|\notag\\
&\leq 2NL^3\sum_{j=i+1}^{N}\|x_j^{[t]}\|.
\label{lem:bound_x:proof_eq2.2}
\end{align}
The second inequality of \eqref{lem:bound_x:proof_eq2} plugs in the bound on $\|\widebar{F}_i(\vtheta^{[t]})-\widebar{F}_i(\vtheta_{1:i-1}^{[t]},\ybf_{i:N}(\vtheta_{1:i-1}^{[t]}))\|$ below, which follows from a similar argument
\begin{align*}
\|\widebar{F}_i(\vtheta^{[t]})-\widebar{F}_i(\vtheta_{1:i-1}^{[t]},\ybf_{i:N}(\vtheta_{1:i-1}^{[t]}))\|&\leq L\|x_i^{[t]}\|+L\sum_{j=i+1}^{N}\|\theta_j^{[t]}-y_j(\vtheta_{1:i-1}^{[t]})\|\notag\\
&\leq 2NL^3\sum_{j=i}^{N}\|x_j^{[t]}\|.
\end{align*}

The term $T_2$ can be bounded leveraging the Lipschitz condition of $y_i$
\begin{align}
\|y_i(\vtheta_{1:i-1}^{[t]})-y_i(\vtheta_{1:i-1}^{[t+1]})\| & \leq L\sum_{j=1}^{i-1}\|\theta_{j}^{[t+1]}-\theta_{j}^{[t]}\|\notag\\
&\leq L\sum_{j=1}^{i-1}\alpha_j^{[t]}\left(\|\Delta f_j^{[t]}\|+NL^2\sum_{k=j}^{N}\|x_k^{[t]}\|\right)\notag\\
&\leq L\sum_{j=1}^{i-1}\alpha_j^{[t]}\|\Delta f_j^{[t]}\|+N^2 L^3\sum_{j=1}^{i-1}\alpha_j^{[t]}\sum_{k=j}^{N}\|x_{k}^{[t]}\|,\label{lem:bound_x:proof_eq4.0}
\end{align}
where the second inequality plugs in the result of Lemma~\ref{lem:Lipschitz}.

Since $\alpha_j^{[t]}\leq\alpha_k^{[t]}$ for any $j\leq k$, ~\eqref{lem:bound_x:proof_eq4.0} implies
\begin{align}
T_2&= \|y_i(\vtheta_{1:i-1}^{[t]})-y_i(\vtheta_{1:i-1}^{[t+1]})\|^2 \notag\\
&\leq 2\left(L\sum_{j=1}^{i-1}\alpha_j^{[t]}\|\Delta f_j^{[t]}\|\right)^2+2\left(N^2 L^3\sum_{j=1}^{N}\alpha_j^{[t]}\|x_{j}^{[t]}\|\right)^2\notag\\
&\leq 2NL^2\sum_{j=1}^{i-1}(\alpha_j^{[t]})^2\|\Delta f_j^{[t]}\|^2+2N^5 L^6\sum_{j=1}^{N}(\alpha_j^{[t]})^2\|x_{j}^{[t]}\|^2.
\label{lem:bound_x:proof_eq4}
\end{align}

Next we treat $T_3$
\begin{align}
T_3 &\leq \alpha_i^{[t]}\|\Delta f_i^{[t]}\|\|x_i^{[t]}-\alpha_i^{[t]}\widebar{F}_i(\vtheta^{[t]})\|\notag\\
&\leq \frac{2\alpha_i^{[t]}}{\delta}\|\Delta f_i^{[t]}\|^2+\frac{\delta\alpha_i^{[t]}}{8}\|x_i^{[t]}-\alpha_i^{[t]}\widebar{F}_i(\vtheta^{[t]})\|^2\notag\\
&\leq \frac{2\alpha_i^{[t]}}{\delta}\|\Delta f_i^{[t]}\|^2+\frac{\delta\alpha_i^{[t]}}{8}\left(\|x_i^{[t]}\|^2-\frac{\delta\alpha_i^{[t]}}{2}\|x_i^{[t]}\|^2+\frac{8N^3 L^6\alpha_i^{[t]}}{\delta}\sum_{j=i+1}^{N}\|x_j^{[t]}\|^2\right)\notag\\
&\leq \frac{2\alpha_i^{[t]}}{\delta}\|\Delta f_i^{[t]}\|^2+\frac{\delta\alpha_i^{[t]}}{8}\|x_i^{[t]}\|^2+\frac{N^3 L^6\alpha_i^{[t]}}{\delta}\sum_{j=i+1}^{N}\|x_j^{[t]}\|^2,
\end{align}
where the third inequality uses \eqref{lem:bound_x:proof_eq2} and the fourth inequality follows from the step size condition $\alpha_i^{[t]}\leq\frac{1}{\delta}$.

To bound $T_4$, we plug in \eqref{lem:bound_x:proof_eq2} and \eqref{lem:bound_x:proof_eq4.0}
\begin{align}
T_4&\leq \|x_i^{[t]}-\alpha_i^{[t]}\widebar{F}_i(\vtheta^{[t]})\|\|y_i(\vtheta_{1:i-1}^{[t]})-y_i(\vtheta_{1:i-1}^{[t+1]})\|\notag\\
&\leq \sqrt{\|x_i^{[t]}\|^2-\frac{\delta\alpha_i^{[t]}}{2}\|x_i^{[t]}\|^2+\frac{8N^3 L^6\alpha_i^{[t]}}{\delta}\sum_{j=i+1}^{N}\|x_j^{[t]}\|^2}\quad\cdot\notag\\
&\hspace{50pt}\left(L\sum_{j=1}^{i-1}\alpha_j^{[t]}\|\Delta f_j^{[t]}\|+N L^3\sum_{j=1}^{i-1}\alpha_j^{[t]}\sum_{k=j}^{N}\|x_k^{[t]}\|\right)\notag\\
&\leq \left(\|x_i^{[t]}\|+\sqrt{\frac{8N^3 L^6\alpha_i^{[t]}}{\delta}}\sum_{j=i+1}^{N}\|x_j^{[t]}\|\right)\cdot\notag\\
&\hspace{50pt}\left(L\sum_{j=1}^{i-1}\alpha_j^{[t]}\|\Delta f_j^{[t]}\|+N^2 L^3\sum_{j=1}^{i-1}\alpha_j^{[t]}\|x_{j}^{[t]}\|+N^2 L^3\alpha_{i-1}^{[t]}\sum_{j=i}^{N}\|x_{j}^{[t]}\|\right)\notag\\
&\leq L\sum_{j=1}^{i-1}\alpha_j^{[t]}\|\Delta f_j^{[t]}\|\|x_i^{[t]}\|+L^4(\alpha_i^{[t]})^{1/2}\alpha_{i-1}^{[t]}\sqrt{\frac{8N^3}{\delta}}\left(\sum_{j=i+1}^{N}\|x_j^{[t]}\|\right)\left(\sum_{j=1}^{i-1}\|\Delta f_j^{[t]}\|\right)\notag\\
&\hspace{20pt}+N^2 L^3\sum_{j=1}^{i-1}\alpha_j^{[t]}\|x_{j}^{[t]}\|\|x_i^{[t]}\|+N^2 L^3\alpha_{i-1}^{[t]} \sum_{j=i}^{N} \|x_{j}^{[t]}\| \|x_i^{[t]}\|\notag\\
&\hspace{20pt}+L^6(\alpha_i^{[t]})^{1/2}\alpha_{i-1}^{[t]}\sqrt{\frac{8N^7}{\delta}}\left(\sum_{j=i+1}^{N}\|x_j^{[t]}\|\right)\left(\sum_{j=1}^{N}\|x_j^{[t]}\|\right)\notag\\
&\leq \frac{NL\alpha_{i-1}^{[t]}}{4}\|x_i^{[t]}\|^2+L\sum_{j=1}^{i-1}\alpha_j^{[t]}\|\Delta f_j^{[t]}\|^2+NL^4(\alpha_i^{[t]})^{1/2}\alpha_{i-1}^{[t]}\sum_{j=i+1}^{N}\|x_j^{[t]}\|^2\notag\\
&\hspace{20pt}+\frac{2N^4 L^4 (\alpha_i^{[t]})^{1/2}\alpha_{i-1}^{[t]}}{\delta}\sum_{j=i+1}^{N}\|\Delta f_j^{[t]}\|^2+\frac{4N^6 L^6\alpha_{i-1}^{[t]}}{\delta}\|x_i^{[t]}\|^2+\frac{\delta}{16N}\sum_{j=1}^{i-1}\alpha_{j}^{[t]}\|x_j^{[t]}\|^2\notag\\
&\hspace{20pt}+\frac{N^3 L^3\alpha_{i-1}^{[t]}}{4}\|x_i^{[t]}\|^2+N^2 L^3\alpha_{i-1}^{[t]}\sum_{j=i}^{N}\|x_j^{[t]}\|^2\notag\\
&\hspace{20pt}+NL^6(\alpha_i^{[t]})^{1/2}\alpha_{i-1}^{[t]}\sqrt{\frac{2N^7}{\delta}}\sum_{j=i+1}^{N}\|x_j^{[t]}\|^2+NL^6(\alpha_i^{[t]})^{1/2}\alpha_{i-1}^{[t]}\sqrt{\frac{2N^7}{\delta}}\sum_{j=1}^{N}\|x_j^{[t]}\|^2\notag\\
&\leq \sum_{j=1}^{i-1}\left(\frac{\delta\alpha_{j}^{[t]}}{16N}+NL^6\sqrt{\frac{2N^7}{\delta}}(\alpha_{i}^{[t]})^2\right)\|x_j^{[t]}\|^2+4N^2L^3\alpha_{i-1}^{[t]}\sum_{j=i}^{N}\|x_j^{[t]}\|^2\notag\\
&\hspace{20pt}+(\frac{N^3L^3\alpha_{i-1}^{[t]}}{2}+\frac{4N^6 L^6\alpha_{i-1}^{[t]}}{\delta})\|x_i^{[t]}\|^2+L\alpha_i^{[t]}\sum_{j=1}^{N}\|\Delta f_j^{[t]}\|^2\notag\\
&\leq \sum_{j=1}^{i-1}\left(\frac{\delta\alpha_{j}^{[t]}}{16N}+NL^6\sqrt{\frac{2N^7}{\delta}}(\alpha_{i}^{[t]})^2\right)\|x_j^{[t]}\|^2+4N^2L^3\alpha_{i-1}^{[t]}\sum_{j=i+1}^{N}\|x_j^{[t]}\|^2\notag\\
&\hspace{20pt}+(\frac{9N^3L^3\alpha_{i-1}^{[t]}}{2}+\frac{4N^6 L^6\alpha_{i-1}^{[t]}}{\delta})\|x_i^{[t]}\|^2+L\alpha_i^{[t]}\sum_{j=1}^{N}\|\Delta f_j^{[t]}\|^2,
\end{align}
where the second inequality follows from the fact that $\sqrt{a_1+a_2+\cdots}\leq\sqrt{a_1}+\sqrt{a_2}+\cdots$ for any positive scalars $a_1,a_2,\cdots$. We have also simplified terms using the step size condition $\alpha_{i}^{[t]}\leq\frac{\delta}{2N^5 L^6}$, $\alpha_{i}^{[t]}\leq\frac{N^2}{L^2}$, and $\alpha_i^{[t]}\leq\frac{\delta^2}{4N^8 L^6}$.

We finally bound $T_5$ with \eqref{lem:bound_x:proof_eq4},
\begin{align}
T_5&\leq \alpha_i^{[t]}\|\Delta f_i^{[t]}\| \|y_i(\vtheta_{1:i-1}^{[t]})-y_i(\vtheta_{1:i-1}^{[t+1]})\|\notag\\
&\leq (\alpha_i^{[t]})^2\|\Delta f_i^{[t]}\|^2+\frac{1}{4} \|y_i(\vtheta_{1:i-1}^{[t]})-y_i(\vtheta_{1:i-1}^{[t+1]})\|^2\notag\\
&\leq (\alpha_i^{[t]})^2\|\Delta f_i^{[t]}\|^2+\frac{NL^2}{2}\sum_{j=1}^{i-1}(\alpha_j^{[t]})^2\|\Delta f_j^{[t]}\|^2+\frac{N^5 L^6}{2}\sum_{j=1}^{N}(\alpha_j^{[t]})^2\|x_{j}^{[t]}\|^2.
\end{align}

Collecting the bounds on $T_1$-$T_5$ and applying them to \eqref{lem:bound_x:proof_eq1}, we have
\allowdisplaybreaks
\begin{align*}
&\|x_i^{t+1}\|^2\notag\\
&\leq \|x_i^{[t]}\|^2-\frac{\delta\alpha_i^{[t]}}{2}\|x_i^{[t]}\|^2+\frac{8N^3 L^6\alpha_i^{[t]}}{\delta}\sum_{j=i+1}^{N}\|x_j^{[t]}\|^2+(\alpha_i^{[t]})^2\|\Delta f_i^{[t]}\|^2\notag\\
&\hspace{20pt}+2NL^2\sum_{j=1}^{i-1}(\alpha_j^{[t]})^2\|\Delta f_j^{[t]}\|^2+2N^5 L^6\sum_{j=1}^{N}(\alpha_j^{[t]})^2\|x_{j}^{[t]}\|^2\notag\\
&\hspace{20pt}+\frac{2\alpha_i^{[t]}}{\delta}\|\Delta f_i^{[t]}\|^2+\frac{\delta\alpha_i^{[t]}}{8}\|x_i^{[t]}\|^2+\frac{N^3 L^6\alpha_i^{[t]}}{\delta}\sum_{j=i+1}^{N}\|x_j^{[t]}\|^2\notag\\
&\hspace{20pt}+\sum_{j=1}^{i-1}\left(\frac{\delta\alpha_{j}^{[t]}}{16N}+NL^6\sqrt{\frac{2N^7}{\delta}}(\alpha_{i}^{[t]})^2\right)\|x_j^{[t]}\|^2+4N^2L^3\alpha_{i-1}^{[t]}\sum_{j=i+1}^{N}\|x_j^{[t]}\|^2\notag\\
&\hspace{20pt}+(\frac{9N^3L^3\alpha_{i-1}^{[t]}}{2}+\frac{4N^6 L^6\alpha_{i-1}^{[t]}}{\delta})\|x_i^{[t]}\|^2+L\alpha_i^{[t]}\sum_{j=1}^{N}\|\Delta f_j^{[t]}\|^2\notag\\
&\hspace{20pt}+(\alpha_i^{[t]})^2\|\Delta f_i^{[t]}\|^2+\frac{NL^2}{2}\sum_{j=1}^{i-1}(\alpha_j^{[t]})^2\|\Delta f_j^{[t]}\|^2+\frac{N^5 L^6}{2}\sum_{j=1}^{N}(\alpha_j^{[t]})^2\|x_{j}^{[t]}\|^2\notag\\
&\leq \|x_i^{[t]}\|^2-\left(\frac{3\delta\alpha_i^{[t]}}{8}-\frac{9N^3L^3\alpha_{i-1}^{[t]}}{2}-\frac{4N^6 L^6\alpha_{i-1}^{[t]}}{\delta}+\frac{5N^5 L^6 (\alpha_i^{[t]})^2}{2}\right)\|x_i^{[t]}\|^2\notag\\
&\hspace{20pt}+\sum_{j=1}^{i-1}\left(\frac{5N^5 L^6 (\alpha_j^{[t]})^2}{2}+\frac{\delta\alpha_{j}^{[t]}}{16N}+NL^6\sqrt{\frac{2N^7}{\delta}}(\alpha_{i}^{[t]})^2\right)\|x_j^{[t]}\|^2\notag\\
&\hspace{20pt}+\sum_{j=i+1}^{N}\left(\frac{9N^3 L^6\alpha_i^{[t]}}{\delta}+\frac{5N^5 L^6 (\alpha_j^{[t]})^2}{2}+4N^2 L^3\alpha_i^{[t]}\right)\|x_j^{[t]}\|^2+(\frac{3}{\delta}+L)\alpha_i^{[t]}\sum_{j=1}^{N}\|\Delta f_j^{[t]}\|^2\notag\\
&\leq \|x_i^{[t]}\|^2-\frac{\delta\alpha_i^{[t]}}{4}\|x_i^{[t]}\|^2+\sum_{j=1}^{i-1}\frac{\delta\alpha_j^{[t]}}{8N}\|x_j^{[t]}\|^2\notag\\
&\hspace{20pt}+\left(\frac{9N^3 L^6}{\delta}+8N^2 L^3\right)\alpha_i^{[t]}\sum_{j=i+1}^{N}\|x_j^{[t]}\|^2+(\frac{3}{\delta}+L)\alpha_i^{[t]}\sum_{j=1}^{N}\|\Delta f_j^{[t]}\|^2,
\end{align*}
where in the second inequality we have combined and simplified terms with the step size condition $\alpha_i^{[t]}\leq\frac{1}{\delta}$ and $\alpha_i^{[t]}\leq\frac{2}{5NL}$, and the third inequality follows from $\frac{\alpha_{i-1}^{[t]}}{\alpha_i^{[t]}}\leq\frac{\delta}{16}(\frac{9N^3 L^3}{2}+\frac{4N^6 L^6}{\delta})^{-1}$, $\alpha_i^{[t]}\leq\frac{\delta}{40N^5 L^6}$, $\frac{(\alpha_{i}^{[t]})^2}{\alpha_1^{[t]}}\leq\frac{\delta^{3/2}}{64N^7}$, $\alpha_i^{[t]}\leq\frac{\delta}{80N^6 L^6}$, and $\frac{(\alpha_{i}^{[t]})^2}{\alpha_1^{[t]}}\leq\frac{8}{5N^3 L^3}$.

\qed

\subsection{Proof of Lemma~\ref{lem:Lipschitz}}

By the definition of $f_i^{[t]}$,
\begin{align}
    \|f_i^{[t]}\|&=\|\Delta f_i^{[t]}+\widebar{F}_i(\vtheta^{[t]})-\widebar{F}_i(\vtheta_{1:i-1}^{[t]},\ybf_{i:N}(\vtheta_{1:i-1}^{[t]}))\|\notag\\
    &\leq\|\Delta f_i^{[t]}\|+\|\widebar{F}_i(\vtheta^{[t]})-\widebar{F}_i(\vtheta_{1:i-1}^{[t]},\ybf_{i:N}(\vtheta_{1:i-1}^{[t]}))\|\notag\\
    &\leq\|\Delta f_i^{[t]}\|+L\sum_{j=i}^{N}\|\theta_j^{[t]}-y_j(\vtheta_{1:i-1}^{[t]})\|.\label{lem:Lipschitz:proof_eq1}
\end{align}

Using the identity \eqref{eq:y_transformation}, we can bound the second term of \eqref{lem:Lipschitz:proof_eq1} recursively
\begin{align}
    &\sum_{j=i}^{N}\|\theta_j^{[t]}-y_j(\vtheta_{1:i-1}^{[t]})\|\notag\\
    &\leq \|x_i^{[t]}\|+\sum_{j=i+1}^{N}\|\theta_j^{[t]}-y_j(\vtheta_{1:i-1}^{[t]})\|\notag\\
    &\leq \|x_i^{[t]}\|+\|\theta_{i+1}^{[t]}-y_{i+1}(\vtheta_{1:i-1}^{[t]}, y_i(\vtheta_{1:i-1}^{[t]}))\|+\sum_{j=i+2}^{N}\|\theta_j^{[t]}-y_j(\vtheta_{1:i-1}^{[t]})\|\notag\\
    &\leq (L+1)\|x_i^{[t]}\|+\|x_{i+1}^{[t+1]}\|+\sum_{j=i+2}^{N}\|\theta_j^{[t]}-y_j(\vtheta_{1:i-1}^{[t]})\|\notag\\
    &\leq (2L+1)\|x_i^{[t]}\|+(L+1)\|x_{i+1}^{[t]}\|+\|x_{i+2}^{[t]}\|+\sum_{j=i+3}^{N}\|\theta_j^{[t]}-y_j(\vtheta_{1:i-1}^{[t]})\|\notag\\
    &\leq NL\sum_{j=i}^{N}\|x_j^{[t]}\|.\label{lem:Lipschitz:proof_eq2}
\end{align}

Plugging \eqref{lem:Lipschitz:proof_eq2} into \eqref{lem:Lipschitz:proof_eq1}, we have
\begin{align}
    \|f_i^{[t]}\|\leq \|\Delta f_i^{[t]}\|+NL^2\sum_{j=i}^{N}\|x_j^{[t]}\|.\label{lem:Lipschitz:proof_eq3}
\end{align}

The bound on $\|\theta_i^{[t+1]}-\theta_i^{[t]}\|$ follows simply from \eqref{lem:Lipschitz:proof_eq3} and the update rule \eqref{alg:update_decision}.

\qed

\subsection{Proof of Lemma~\ref{lem:conv_x_standardSA}}
\label{sec:proof_lem:conv_x_standardSA}

By the update rule \eqref{alg:update_standardSA}, we have
\begin{align*}
x_i^{[t+1]}&=\theta_{i}^{[t+1]}-y_i(\vtheta_{1:i-1}^{[t+1]})\notag\\
&=\theta_{i}^{[t]}-\alpha_i^{[t]}F_i(\vtheta^{[t]},X^{[t]})-y_i(\vtheta_{1:i-1}^{[t+1]})\notag\\
&=x_i^{[t]}-\alpha_i^{[t]}F_i(\vtheta^{[t]},X^{[t]})+\left(y_i(\vtheta_{1:i-1}^{[t]})-y_i(\vtheta_{1:i-1}^{[t+1]})\right)\notag\\
&=x_i^{[t]}-\alpha_i^{[t]} \widebar{F}_i(\vtheta^{[t]})-\alpha_i^{[t]}\big(F_i(\vtheta^{[t]},X^{[t]})-\widebar{F}_i(\vtheta^{[t]})\big)+\big(y_i(\vtheta_{1:i-1}^{[t]})-y_i(\vtheta_{1:i-1}^{[t+1]})\big).
\end{align*}

Taking the norm,
\begin{align}
&\|x_i^{[t+1]}\|^2\notag\\
&=\underbrace{\|x_i^{[t]}-\alpha_i^{[t]}\widebar{F}_i(\vtheta^{[t]})\|^2}_{A_1}+\underbrace{(\alpha_i^{[t]})^2\|F_i(\vtheta^{[t]},X^{[t]})-\widebar{F}_i(\vtheta^{[t]})\|^2}_{A_2}\notag\\
&\hspace{20pt}+\underbrace{\|y_i(\vtheta_{1:i-1}^{[t]})-y_i(\vtheta_{1:i-1}^{[t+1]})\|^2}_{A_3}+\underbrace{\alpha_i^{[t]}\langle x_i^{[t]}-\alpha_i^{[t]}\widebar{F}_i(\vtheta^{[t]}),F_i(\vtheta^{[t]},X^{[t]})-\widebar{F}_i(\vtheta^{[t]})\rangle}_{A_4}\notag\\
&\hspace{20pt}+\underbrace{\langle x_i^{[t]}-\alpha_i^{[t]}\widebar{F}_i(\vtheta^{[t]}),y_i(\vtheta_{1:i-1}^{[t]})-y_i(\vtheta_{1:i-1}^{[t+1]})\rangle}_{A_5}\notag\\
&\hspace{20pt}+\underbrace{\alpha_i^{[t]}\langle F_i(\vtheta^{[t]},X^{[t]})-\widebar{F}_i(\vtheta^{[t]}), y_i(\vtheta_{1:i-1}^{[t]})-y_i(\vtheta_{1:i-1}^{[t+1]})}_{A_6}\rangle.
\label{lem:conv_x_standardSA_proof_eq1}
\end{align}

A bound on $A_1$ can be obtained using an argument identical to \eqref{lem:bound_x:proof_eq2}, provided that the step sizes satisfy $\alpha_i^{[t]}\leq \frac{\delta}{8N^3 L^6}$ and $\alpha_i^{[t]}\leq \frac{1}{\delta}$
\begin{align}
A_1&=\|x_i^{[t]}\|^2-2\alpha_i^{[t]}\langle x_i^{[t]},\widebar{F}_i(\vtheta^{[t]})\rangle+(\alpha_i^{[t]})^2\|\widebar{F}_i(\vtheta^{[t]})\|^2\notag\\
&\leq \|x_i^{[t]}\|^2-\frac{\delta\alpha_i^{[t]}}{2}\|x_i^{[t]}\|^2+\frac{8N^3 L^6\alpha_i^{[t]}}{\delta}\sum_{j=i+1}^{N}\|x_j^{[t]}\|^2.
\label{lem:conv_x_standardSA_proof_eq2}
\end{align}


The second term is bounded in expectation due to Assumption~\ref{assump:bounded}
\begin{align}
A_2=(\alpha_i^{[t]})^2\|F_i(\vtheta^{[t]},X^{[t]})-\widebar{F}_i(\vtheta^{[t]})\|^2\leq 4D^2(\alpha_i^{[t]})^2.\label{lem:conv_x_standardSA_proof_eq3}
\end{align}

Using the Lipschitz continuity of $y_i$, we have
\begin{align}
\|y_i(\vtheta_{1:i-1}^{[t]})-y_i(\vtheta_{1:i-1}^{[t+1]})\|& \leq L\sum_{j=1}^{i-1}\|\theta_{j}^{[t+1]}-\theta_{j}^{[t]}\|\notag\\
&= L\sum_{j=1}^{i-1}\alpha_j^{[t]}\|F_i(\vtheta^{[t]},X^{[t]})\|\notag\\
&\leq NLD\alpha_{i-1}^{[t]},\label{lem:conv_x_standardSA_proof_eq4}
\end{align}
which implies
\begin{align}
A_3&= \|y_i(\vtheta_{1:i-1}^{[t]})-y_i(\vtheta_{1:i-1}^{[t+1]})\|^2 \leq N^2 L^2 D^2 (\alpha_{i-1}^{[t]})^2.
\label{lem:conv_x_standardSA_proof_eq5}
\end{align}

The term $A_4$ is obviously zero in expectation as
\begin{align*}
\mathbb{E}[F_i(\vtheta^{[t]},X^{[t]})-\widebar{F}_i(\vtheta^{[t]})]=0.
\end{align*}

The bound on $A_5$ follows from the Cauchy-Schwarz inequality
\begin{align}
A_5&\leq\|x_i^{[t]}\|\|y_i(\vtheta_{1:i-1}^{[t]})-y_i(\vtheta_{1:i-1}^{[t+1]})\|+\alpha_i^{[t]}\|\widebar{F}_i(\vtheta^{[t]})\|\|y_i(\vtheta_{1:i-1}^{[t]})-y_i(\vtheta_{1:i-1}^{[t+1]})\|\notag\\
&\leq \frac{\delta\alpha_i^{[t]}}{4}\|x_i^{[t]}\|^2+\frac{1}{\delta\alpha_i^{[t]}}\|y_i(\vtheta_{1:i-1}^{[t]})-y_i(\vtheta_{1:i-1}^{[t+1]})\|^2+D\alpha_i^{[t]}\|y_i(\vtheta_{1:i-1}^{[t]})-y_i(\vtheta_{1:i-1}^{[t+1]})\|\notag\\
&\leq \frac{\delta\alpha_i^{[t]}}{4}\|x_i^{[t]}\|^2+\frac{(NLD\alpha_{i-1}^{[t]})^2}{\delta\alpha_i^{[t]}}+D\alpha_i^{[t]}\cdot NLD\alpha_{i-1}^{[t]}\notag\\
&\leq \frac{\delta\alpha_i^{[t]}}{4}\|x_i^{[t]}\|^2+\frac{N^2 L^2 D^2 (\alpha_{i-1}^{[t]})^2}{\delta\alpha_i^{[t]}}+ NLD^2\alpha_i^{[t]}\alpha_{i-1}^{[t]}
\label{lem:conv_x_standardSA_proof_eq7}
\end{align}
where in the third equality we plug in \eqref{lem:conv_x_standardSA_proof_eq4}.

Finally, we bound $A_6$ using a similar argument
\begin{align}
A_6&\leq \alpha_i^{[t]}\|F_i(\vtheta^{[t]},X^{[t]})-\widebar{F}_i(\vtheta^{[t]})\| \|y_i(\vtheta_{1:i-1}^{[t]})-y_i(\vtheta_{1:i-1}^{[t+1]})\|\notag\\
&\leq \alpha_i^{[t]}\cdot 2D\cdot NLD\alpha_{i-1}^{[t]}\notag\\
&\leq  2 NLD^2\alpha_{i-1}^{[t]}\alpha_i^{[t]}.
\label{lem:conv_x_standardSA_proof_eq8}
\end{align}

Collecting the bounds in \eqref{lem:conv_x_standardSA_proof_eq2}-\eqref{lem:conv_x_standardSA_proof_eq8} and putting them into \eqref{lem:conv_x_standardSA_proof_eq1}, we get
\begin{align*}
&\mathbb{E}[\|x_i^{[t+1]}\|^2]\\
&\leq \mathbb{E}[\|x_i^{[t]}\|^2-\frac{\delta\alpha_i^{[t]}}{2}\|x_i^{[t]}\|^2+\frac{8N^3 L^6\alpha_i^{[t]}}{\delta}\sum_{j=i+1}^{N}\|x_j^{[t]}\|^2+4D^2(\alpha_i^{[t]})^2+N^2 L^2 D^2 (\alpha_{i-1}^{[t]})^2\notag\\
&\hspace{20pt}+\frac{\delta\alpha_i^{[t]}}{4}\|x_i^{[t]}\|^2+\frac{N^2 L^2 D^2 (\alpha_{i-1}^{[t]})^2}{\delta\alpha_i^{[t]}}+ NLD^2\alpha_i^{[t]}\alpha_{i-1}^{[t]}+2 NLD^2\alpha_{i-1}^{[t]}\alpha_i^{[t]}]\notag\\
&\leq \mathbb{E}[\|x_i^{[t]}\|^2]-\frac{\delta\alpha_i^{[t]}}{4}\mathbb{E}[\|x_i^{[t]}\|^2]+\frac{8N^3 L^6\alpha_i^{[t]}}{\delta}\sum_{j=i+1}^{N}\mathbb{E}[\|x_j^{[t]}\|^2]\notag\\
&\hspace{20pt}+8N^2 L^2 D^2 (\alpha_{i}^{[t]})^2+\frac{N^2 L^2 D^2 (\alpha_{i-1}^{[t]})^2}{\delta\alpha_i^{[t]}}.
\end{align*}

\qed

\subsection{Proof of Lemma~\ref{lem:Markovian_crossterm}}

This lemma essentially bounds the distance between samples from a time-varying Markov chain and the stationary distribution. Techniques for proving this lemma has been well-studied in the literature. Hence we omit the proof here but note that it is very similar to the proof of Lemma 1 in \citet{zeng2024two}.

\end{document}